\documentclass{article}
\usepackage{amsmath,amsthm,amsfonts,amssymb,amscd,latexsym, fullpage,textcomp,amscd}
\date{ }
\begin{document}

\theoremstyle{plain}
\newtheorem{thm}{\bf Theorem}[section]
\newtheorem{lem}[thm]{\bf Lemma}
\newtheorem{cor}[thm]{\bf Corollary}
\newtheorem{prop}[thm]{\bf Proposition}
\theoremstyle{remark}
\newtheorem{Case}{\bf Case}
\newtheorem{rem}[thm]{\bf Remark}
\newtheorem{claim}[thm]{\bf Claim}

\theoremstyle{definition}
\newtheorem{Def}{\bf Definition}
\newtheorem*{pf}{\bf Proof}
\newtheorem{Conj}{\bf Conjecture}

\newcommand{\nc}{\newcommand}

\newcommand{\BB}{{\mathcal B}} 
\newcommand{\CC}{{\mathcal C}}
\newcommand{\DD}{{\mathcal D}}
\newcommand{\EE}{{\mathcal E}}
\newcommand{\FF}{{\mathcal F}}
\newcommand{\GG}{{\mathcal G}}
\newcommand{\HH}{{\mathcal H}}
\newcommand{\II}{{\mathcal I}}
\newcommand{\JJ}{{\mathcal J}}
\newcommand{\KK}{{\mathcal K}}
\newcommand{\LL}{{\mathcal L}}
\newcommand{\MM}{{\mathcal M}}
\newcommand{\NN}{{\mathcal N}}
\newcommand{\OO}{{\mathcal O}}
\newcommand{\PP}{{\mathcal P}}
\newcommand{\QQ}{{\mathcal Q}}
\newcommand{\RR}{{\mathcal R}}
\newcommand{\TT}{{\mathcal T}}
\newcommand{\UU}{{\mathcal U}}
\newcommand{\VV}{{\mathcal V}}
\newcommand{\WW}{{\mathcal W}}
\newcommand{\ZZ}{{\mathcal Z}}
\newcommand{\XX}{{\mathcal X}}
\newcommand{\YY}{{\mathcal Y}}
\nc{\bba}{{\mathbb A}}
\nc{\bbb}{{\mathbb B}}
\nc{\bbc}{{\mathbb C}}
\nc{\bbd}{{\mathbb D}}
\nc{\bbe}{{\mathbb E}}
\nc{\bbf}{{\mathbb F}}
\nc{\bbg}{{\mathbb G}}
\nc{\bbh}{{\mathbb H}}
\nc{\bbi}{{\mathbb I}}
\nc{\bbj}{{\mathbb J}}
\nc{\bbk}{{\mathbb K}}
\nc{\bbl}{{\mathbb L}}
\nc{\bbm}{{\mathbb M}}
\nc{\bbo}{{\mathbb O}}
\nc{\bbp}{{\mathbb P}}
\nc{\bbq}{{\mathbb Q}}
\nc{\bbr}{{\mathbb R}}
\nc{\bbs}{{\mathbb S}}
\nc{\bb}{{\mathbb T}}
\nc{\bbu}{{\mathbb U}}
\nc{\bbv}{{\mathbb V}}
\nc{\bbw}{{\mathbb W}}
\nc{\bbx}{{\mathbb X}}
\nc{\bby}{{\mathbb Y}}
\nc{\bbz}{{\mathbb Z}}
\nc{\fA}{{\mathfrak A}}
\nc{\fB}{{\mathfrak B}}
\nc{\fC}{{\mathfrak C}}
\nc{\fD}{{\mathfrak D}}
\nc{\fE}{{\mathfrak E}}
\nc{\fF}{{\mathfrak F}}
\nc{\fG}{{\mathfrak G}}
\nc{\fH}{{\mathfrak H}}
\nc{\fI}{{\mathfrak I}}
\nc{\fJ}{{\mathfrak J}}
\nc{\fK}{{\mathfrak K}}
\nc{\fL}{{\mathfrak L}}
\nc{\fM}{{\mathfrak M}}
\nc{\fN}{{\mathfrak N}}
\nc{\fO}{{\mathfrak O}}
\nc{\fP}{{\mathfrak P}}
\nc{\fQ}{{\mathfrak Q}}
\nc{\fR}{{\mathfrak R}}
\nc{\fS}{{\mathfrak S}}
\nc{\fT}{{\mathfrak T}}
\nc{\fU}{{\mathfrak U}}
\nc{\fV}{{\mathfrak V}}
\nc{\fW}{{\mathfrak W}}
\nc{\fZ}{{\mathfrak Z}}
\nc{\fX}{{\mathfrak X}}
\nc{\fY}{{\mathfrak Y}}
\nc{\fa}{{\mathfrak a}}
\nc{\fb}{{\mathfrak b}}
\nc{\fc}{{\mathfrak c}}
\nc{\fd}{{\mathfrak d}}
\nc{\fe}{{\mathfrak e}}
\nc{\ff}{{\mathfrak f}}
\nc{\fh}{{\mathfrak h}}
\nc{\fj}{{\mathfrak j}}
\nc{\fk}{{\mathfrak k}}
\nc{\fl}{{\mathfrak{l}}}
\nc{\fm}{{\mathfrak m}}
\nc{\fn}{{\mathfrak n}}
\nc{\fo}{{\mathfrak o}}
\nc{\fp}{{\mathfrak p}}
\nc{\fq}{{\mathfrak q}}
\nc{\fr}{{\mathfrak r}}
\nc{\fs}{{\mathfrak s}}
\nc{\ft}{{\mathfrak t}}
\nc{\fu}{{\mathfrak u}}
\nc{\fv}{{\mathfrak v}}
\nc{\fw}{{\mathfrak w}}
\nc{\fz}{{\mathfrak z}}
\nc{\fx}{{\mathfrak x}}
\nc{\fy}{{\mathfrak y}}

\nc{\al}{{\alpha }}
\nc{\be}{{\beta }}
\nc{\ga}{{\gamma }}
\nc{\de}{{\delta }}
\nc{\vap}{{\tepsilon }}

\nc{\ze}{{\zeta }}
\nc{\et}{{\eta }}
\nc{\vth}{{\vartheta }}

\nc{\io}{{\iota }}
\nc{\ka}{{\kappa }}
\nc{\la}{{\lambda }}
\nc{\vpi}{{     \varpi          }}
\nc{\vrho}{{    \varrho         }}
\nc{\si}{{      \sigma          }}
\nc{\ups}{{     \upsilon        }}
\nc{\vphi}{{    \varphi         }}
\nc{\om}{{      \omega          }}

\nc{\Ga}{{\Gamma }}
\nc{\De}{{\Delta }}
\nc{\nab}{{\nabla}}
\nc{\Th}{{\Theta }}
\nc{\La}{{\Lambda }}
\nc{\Si}{{\Sigma }}
\nc{\Ups}{{\Upsilon }}
\nc{\Om}{{\Omega }}

\nc{\zz}{{\mathbb Z}}
\newcommand{\N}{{\mathbb N}}
\newcommand{\etat}{\tilde{\eta}}
\newcommand{\sln}{\mathfrak{sl} _n }
\newcommand{\slnr}{{\mathfrak {sl} _n (\mathbb R)}}
\newcommand{\sun}{\mathfrak{su} _n}
\newcommand{\cc}{{\mathbb C}}
\newcommand{\rr}{{\mathbb R}}
\newcommand{\ac}{{\check{\alpha}}}
\newcommand{\orb}{{\mathcal O}}
\newcommand{\gcs}{{\mathcal{GC} _{G/K}}}
\newcommand{\cds}{{\mathbb{C} \mathcal{D} _{G/K}}}
\newcommand{\ds}{{\mathcal{D} _{G/K}}}
\newcommand{\ggcs}{{\mathcal{GC} _{G/K} ^G}}
\newcommand{\gcds}{{\mathbb{C} \mathcal{D} _{G/K} ^G} }
\newcommand{\gds}{{\mathcal{D} _{G/K} ^G}}
\newcommand{\bc}{{\check{\beta}}}
\newcommand{\tep}{{\varepsilon_\sharp}}
\newcommand{\ep}{{\varepsilon_\sharp}}
\newcommand{\vep}{\varepsilon}
\newcommand{\epp}{T_{\overline{\varepsilon}}}
\newcommand{\ctsd}{\mathfrak{g} _\mathbb{C} ^*}
\newcommand{\gts}{$\mathfrak{g}\oplus\mathfrak{g}^{\ast}$}
\newcommand{\inv}{^{-1}}
\newcommand{\fg}{{\mathfrak g}}
\newcommand{\noi}{\noindent}
\newcommand{\lra}{{\longrightarrow}}
\newcommand{\pbd}{f^\star \mathcal{D}}
\newcommand{\cep}{{\overline{\varepsilon}}}
\newcommand{\cdt}{\tilde{\mathcal{D}}}
\newcommand{\cts}{{\mathfrak{g} _{\mathbb{C}}}}
\newcommand{\csts}{\mathfrak{k} _{\mathbb{C}}}
\newcommand{\cgts}{{\mathfrak{g} _{\mathbb C} \oplus \mathfrak{g} _\mathbb{C} ^* } }

\newcommand{\pionep}{{ (\pi_ 1)_\sharp}}
\newcommand{\pitwop}{{(\pi _2 )_\sharp}}
\newcommand{\omonep}{{(\omega _1)_\sharp}}
\newcommand{\omtwop}{{(\omega_2)_\sharp}}
\nc{\st}{{\; | \; }}

\nc{\trm}{\textreferencemark}
\nc{\tih}{{\tilde \HH}}
\nc{\vb}{\text{vector bundle}}
\nc{\vbm}{{\VV ect \BB und _M}}
\nc{\vs}{{Vect_{fd}}}
\nc{\cxn}{\text{connection}}
\nc{\bla}{{\mathfrak g ^* \rtimes \mathfrak g}}
\nc{\cbla}{{\mathfrak g _\cc ^* \rtimes \mathfrak g _\cc}}
\nc{\gc}{\text{generalized complex }}
\nc{\gcstr}{\text{generalized complex structure }}
\nc{\gcstrs}{\text{generalized complex structures }}

\title{Two Categories of Dirac Manifolds}
\maketitle
\begin{center} 
\author{Brett Milburn}
\end{center}

\begin{abstract} \noi Generalized complex geometry \cite{hit} and, more generally, Dirac geometry~\cite{cou}, \cite{cow}, unify several familiar geometric structures into one uniform viewpoint.  We introduce two new notions of morphisms between manifolds equipped with Dirac structures, giving two different Dirac categories.  The first generalizes holomorphic and Poisson maps, while the second dual notion generalizes both symplectic and holomorphic maps.  As an application, we consider Dirac groups (i.e. Lie groups with Dirac structure such that group multiplication is a Dirac map).  We explain the conditions under which a group with Dirac structure is a Dirac group.  More precisely, we explain the data and conditions for a Dirac group.  Dirac groups turn out to be a generalization of Poisson groups. 
\end{abstract}






\tableofcontents

\section{Introduction} 
 
\noi Complex Dirac structures on manifolds
generalize several basic
concepts of differential geometry:\
Poisson, presymplectic, and complex structures as well as 
integrable distributions.  Dirac structures are originally due to Courant and Weinstein \cite{cow}, \cite{cou}.  Formally, this is a generalization of
integrable distributions 
where the tangent bundle $TM$ is replaced
by $  \VV _M=TM\oplus T^*M$ and the Lie algebroid structure 
on $TM$ is replaced by
the Courant algebroid structure which is a bracket operation $[-,-]$ on
$ \VV _M$.  Hitchin's notion of a generalized complex structure~\cite{hit} is a specific type of Dirac structure on which we will focus much of our attention here.  Generalized complex geometry contains both complex and symplectic geometry and has found a number of applications in mathematics and physics 
~\cite{kap},~\cite{hit},~\cite{gua}.  \\

\noi Although Dirac geometry is a common viewpoint from which to consider these more traditional structures, there has not yet been a definition of morphisms between manfiolds with Dirac structure which has been widely used or accepted.  Several possible approaches have been put forth and are described in section~\ref{comp}.   In this paper we define two categories of Dirac manifolds (i.e. manifolds with complex Dirac structures).  The first notion of maps we call \emph{Dirac
  maps}, and the corresponding category is seen to 
contain the
categories of Poisson and complex manifolds as full subcategories. We call this category simply the \emph{Dirac category}.  Actually, the Dirac category seems to be exactly a unification or common setting for Poisson and complex manifolds.  With Poisson and complex maps at two opposite ends of the spectrum, the general case is an interpolation between or mixture of the two.\\

\noi Our only application of Dirac maps at present is the study of groups in the category of Dirac manifolds.  These \emph{Dirac groups} are a combination of Poisson groups and complex groups, though it is useful to think of them as a generalization of Poisson groups.  A particularly interesting class of Dirac groups are generalized complex groups.  These turn out to be equivalent to holomorphic Poisson groups, so they provide a (non-holomorphic) point of view on quantization of holomorphic Poisson groups.\\

\noi The
second notion, \emph{dual-Dirac maps}, defines a \emph{dual-Dirac category}
which contains presymplectic and complex manifolds as full
subcategories. Generalized complex structures can be viewed as an interpolation between symplectic and complex structures, especially in the hyperk\"ahler setting~\cite{gua}.  The dual-dirac maps allow for a uniform point of view on symplectic and complex maps so that in any setting involving a variation between the two extremes, it makes sense to talk about maps in a consistent way.   The notion of  dual-Dirac maps has an additional property that
it is stable under
B-transforms.  \\

\noi Dirac maps and dual-Dirac maps provide two new structures of a category on Hitchin's generalized complex manifolds, i.e. two reasonable notions of generalized complex maps.  We establish conditions for which multiplication in a Lie group is a dual-Dirac map and partially classify group objects in the dual-Dirac category.  \\

\noi Finally, we generalize further by also considering categories of Dirac manifolds for which Dirac structures lie in arbitrary exact Courant algebroids (i.e. Courant algebroids $\EE$ on $M$ which are extensions $0 \lra T^*M \lra \EE \lra TM \lra 0$ of the tangent Lie algebroid by the cotangent bundle).  This generality is crucial for expected applications to Representation Theory of affine Lie algebras.  We mostly consider the dual-Dirac category, which extends readily and naturally to the case of arbitrary exact Courant algebroids because of stability under B-transforms.  We also find an extension of the Dirac category to exact Courant algebroids.  However, this extension is less natural since it requires keeping track of additional data.  The problem of a natural theory of Dirac maps for exact Courant algebroids seems to have real content.   It may be related to the notion of \emph{vertex algebroids}, which are a quantization of the notion of Courant algebroids.  \\



\section{Recollections on Dirac and Generalized Complex Geometry}\label{expo} 

\noi This section is a brief expostion on the essential ideas in generalized complex geometry.  Here we introduce the basic definitions and notational conventions used in
this paper.  For a systematic development of generalized complex
structures as well as some of their applications, we refer the reader to ~\cite{gua}. \\

\noi For a manifold $M$, generalized geometry is concerned with the bundle
$\mathcal V _M := TM \oplus
T^*M$.  There is a natural bilinear form on $\mathcal V _M$, given by
the obvious pairing $\langle X + \xi , Y + \eta \rangle = X(\eta) + Y (\xi)$ for sections $X,Y$ of $TM$ and $\xi $, $\eta$ of $T^*M$.  Furthermore,
$\mathcal V _M$ is equipped with the \emph{Courant bracket} defined by 
\[ [X + \xi , Y + \eta ] = [X,Y] + \io _X d \eta +\frac{1}{2} d( \io _X \eta) - \io _Y d \xi - \frac{1}{2} d(\io _Y \xi),\]
where $\io$ denotes contraction in the first variable ($ \io _x \phi = \phi (x,-,...)$).  
The Courant bracket $[ \; , \; ]$ and the bilinear form $\langle \; , \; \rangle$
extend $\cc$-bilinearly to $(\mathcal V _M )_\cc = \mathcal V _M
\otimes \cc$.  


\begin{Def} A \emph{generalized almost complex structure} on $M$ is a map $\JJ : \VV _M \lra \VV _M$ such that $\JJ$ is orthogonal with respect to the inner product $\langle \; , \; \rangle$ and $\JJ ^2 = -1$.  Just as with complex structures, one may consider the $i$-eigenbundle, $D$, of $\JJ $ in $(\VV _M)_\cc$.  The Courant bracket defines an integrability condition ($[D , D] \subset D$) for $\JJ$ to be called a \emph{ generalized complex structure}.  This follows the analogy with almost complex stuctures; an almost complex structure is a complex structure precisely when its $i$-eigenbundle is integrable with respect to the Lie bracket. 
\end{Def} 

\noi The two canonical examples of generalized complex structures come from complex and symplectic structures.  Since $\VV _M = TM \oplus T^*M$, we can express any map $\VV _M \lra \VV _M$ as a block matrix in terms of this decomposition, and we will follow this convention throughout the text.  If $J$ is a complex structure, 

\[
\begin{bmatrix}
J & 0 \\
0 & -J^*
\end{bmatrix} 
\]

\noi is a generalized complex structure.  The $i $-eigenbundle of $\JJ$ is $E \oplus Ann(E)$, where $E$ is the $i$-eigenbundle of $J$, and $Ann(E)$ is the annihilator of $E$ in $T^*M$.  \\

\noi For a symplectic structure $\om$, we get a generalized complex structure 
\[
\begin{bmatrix}
0 & -\om _\sharp \inv \\
\om _\sharp & 0
\end{bmatrix} 
\]
\noi where $\om _\sharp (x) := \om (x , -)$.  The $i$-eigenbundle is the graph of $i\om _\sharp $ in $(\VV _M)_\cc$.  The fact that the symplectic form $\om$ is closed implies that this generalized almost complex structure is integrable, hence a generalized complex structure.   \\

\noi The $i$-eigenbundle $D$ of a generalized complex structure $\JJ$ turns out to be an integrable maximal isotropic subbundle of $(\VV _M)_\cc$, also known as a \emph{complex Dirac structure}.  
Thus, the study of generalized geometry now lies in the framework of Dirac structures.   With this in mind, we recall the following working definitions for our paper.  

\begin{Def} For any manifold, $M$, 
\begin{enumerate}
\item  A \emph{real almost Dirac structure} on $M$ is a maximal isotropic
  subbundle $D$ of  $\mathcal V _M$. A real almost Dirac structure is called a \emph{real Dirac structure} if it is
  integrable with respect to the Courant bracket.  Similarly, a \emph{complex almost Dirac structure} is a maximal isotropic
  subbundle $D \subset (\mathcal V _M )_\cc $, and a \emph{complex Dirac
  structure} is an integrable complex almost Dirac structure.   
\item A complex Dirac structure $D$ is said to be of \emph{constant rank} if the projection map $pr : D \lra TM$ is of constant rank.  
\end{enumerate}
\end{Def}

\noi A generalized (almost) complex structure $\JJ$ is equivalent to is a complex (almost) Dirac structure $D$
  such that $D \cap \overline D = 0$.  Note that the integrability is a closed condition and that being generalized complex is an open condition.  Henceforth we will think of generalized complex structures as complex Dirac structures.  \\

\noi Since the complexification of any Dirac structure is a complex Dirac
structure, both generalized complex structures and Dirac structures
are complex Dirac structures.  Thus, the set of complex Dirac structures contains real Dirac structures and generalized complex structures.  Henceforth (almost) \emph{Dirac structure} will always mean complex (almost) Dirac structure, and we will specify whether it is also real Dirac (i.e. if $\overline  {\mathcal D } = \mathcal D$) if there is any ambiguity. \\

\noi Most of the complex Dirac structures
considered in this paper will be of constant rank. 
It is checked in ~\cite{gua} that any complex Dirac structure of constant rank is of the form 
\[ L(E , \vep ) := \{ X + \xi \in (\mathcal V _M)_\cc \; | \; X \in E \; and \;  \iota _X \vep = \xi _{|E} \}, 
\]
\noi where $E$ is a subbundle of $TM$ and $\vep \in \Ga (M , \wedge ^2 E^*)$.  Complex and symplectic structures, for example, are of this form.  For a subbundle $E$ of $TM$, we define the differential \\
$d_E : \Ga (M, \wedge ^2 E^*) \lra \Ga (M, \wedge ^3 E^*)$ by the following forumula.  For sections $X,Y,Z$ of $E$ and  $\vep \in \Ga (M , \wedge ^2 E^*)$,  
\[d_E \vep (X,Y,Z) = \vep (X,[Y,Z]) + \vep (Y,[Z,X]) + \vep (Z,[X,Y]) + X\vep (Y,Z) - Y\vep (X,Z) + Z\vep (X,Y). \]
\noi In other words, $d_E \vep $ is the restriction to $\wedge ^3 E $ of the ordinary De Rham differential of any extension $\tilde \vep \in \wedge ^2 T^*M $ of $\vep$.  Gualtieri \cite{gua} proves the following useful lemma.

\begin{lem}\label{21sept1} A complex almost  Dirac structure of constant rank $L(E ,
  \vep )$ is integrable if and only if $E$ is integrable and $d_E \vep
  =0$.  
\end{lem}  

\begin{rem} Dirac groups (\S \ref{13Dec4}) provide interesting examples of Dirac structures which are not of constant rank.  In particular, generalized complex groups provide examples of generalized complex structures not of constant rank.   
\end{rem}  

\begin{lem}We recall from~\cite{gua} that the following procedures create Dirac structures 
from geometric structures on $M$.  The first three are real, and (4) and (5) are special cases of generalized complex structures (as we have seen).
\begin{enumerate}
\item 
To an integrable distribution $\mathcal{D}\subset TM$, assign
$[\mathcal{D}\oplus Ann(\mathcal{D})]_\cc$.
\item 
To a Poisson structure $\pi \in \Ga (M, \wedge ^2 TM)$, assign $L(\pi , T^*M)$.
\item 
To a presymplectic structure $\om \in \Om ^2 (M)$, assign $L(TM , \om)$.

\item 
To a complex structure
$J$, assign $T^{(1,0)}M \oplus T^{* ,(0,1)}M$, where $T^{(1,0)} M$ and $T^{(0,1)}M$ denote the holomorphic and antiholomorphic tangent bundles respectively with respect to $J$. 
\item To a symplectic structure $\om \in \Om ^2 (M)$, assign $L(TM_\cc , i \om)$.  
\end{enumerate}
\end{lem}

\begin{rem} A symplectic form $\omega$ on $M$  determines a complex Dirac structure
in one of two ways: $L(TM, \omega) $ and $L(T_\cc M, -i\omega)$.  The former
is a real Dirac structure, and the latter is a generalized complex
structure. 
\end{rem}

\noi This way of representing a Dirac structure $L$ as $L(E,\vep)$ turns out to be extremely useful for our purposes.  
For any vector bundle $E$ and $\vep \in \bigwedge ^2 E^*$, the convention used in this paper is for $\varepsilon _\sharp$ to denote the map $E \lra E^*$ determined by $\vep$.  That is, for $X,Y \in E$, $(\vep _\sharp X)(Y) = \vep (X,Y)$.  \\

\noi For a Dirac structure $D$, if the projection $pr : D \lra T^*M$ has constant rank, then there is some subbundle $U \subset T^*M$ and some $\pi \in \Ga (M, \bigwedge ^2 U^*)$ such that $D$ is of the form 
\[ L(\pi , U) := \{ X + \xi \st X_{|U} = \io _\xi \pi \}.  \] 
If $U = T^*M$, then for $\pi \in \Ga (M , \wedge ^2 TM)$, $L(\pi , T^*M)$ is a Dirac structure if and only if $\pi$ is a Poisson bi-vector ~\cite{gua}, ~\cite{vai}.  Now presymplectic structures, complex structures, and Poisson structures can all be considered Dirac structures.  \\

\noi We recall the notions of pullback and pushforward of linear Dirac structures ~\cite{gua}. For a map $F: V \lra W$ of vector spaces and a subspace $D \subset V \oplus V^*$, define
\[F_\star D = \{ FX + \xi \in W \oplus W^* \st X + F^* \xi \in D \} \]
and for a subspace $D \subset W \oplus W^*$, define 
\[ F^\star D = \{ X + F^* \xi \in V \oplus V \st FX + \xi \in D \} = (F^*)_\star D. \]

\noi Now let $f : M \lra N$ be any map of
manifolds.  For a Dirac structure, $D$, on $N$, the pullback $f^\star D$  is
defined pointwise by $(f^\star D)_p = (df_p)^\star D_{f(p)} $. It is
not necessarily itself a Dirac structure.  \\

\subsection{Twisted Courant Bracket and Automorphisms}

\noi In addition to the standard Courant bracket on $\VV _M$, \v Severa and Weinstein noticed a twisted Courant bracket $[ \; , \; ] _H$ for each closed 3-form $H$ on $M$, defined as
\[ [X + \xi , Y + \eta  ] _H = [X + \xi , Y + \eta  ] + H(X,Y,-) 
\]
and $\VV _{M,H} := (TM \oplus T^*M , [ \; , \; ] _H , \langle \; , \; \rangle )$ so that $\VV _M$ is just $\VV _{M,0}$.  For any 2-form $B$ on $M$, there is an automorphism of the vector bundle $TM \oplus T^*M$, 

\[
\begin{bmatrix}
1 & 0 \\
B_\sharp & 1
\end{bmatrix}, 
\]
denoted $e^B$.  
Indeed, $e^B$ is an isomorphism $\VV _{M,H+dB} \lra \VV _{M, H}$.  In other words, $e^B$ is orthogonal with respect to $ \langle \; , \; \rangle $, and $[e^B u , e^B v]_H = e^B [u,v]_{H +dB}$ for all sections $u,v$ of $TM \oplus T^*M$.  When $B$ is closed, then $e^B$ is an automorphism of $\VV _{M,H}$.  This is what we call a \emph{B-transform}, i.e. an automorphism of $\VV _{M,H}$ of the form $e^B$ for a closed 2-form $B$.   In fact, the automorphism group of $\VV _{M,H}$ is the semidirect product of the group of diffeomorphisms $M \lra M$ and closed 2-forms $Z^2(M)$ \cite{gua}.  B-transforms are thought of as the symmetries of the Courant bracket. \\

\noi An \emph{H-twisted Dirac structure} $D \subset \VV _{M,H}$ is simply a maximal isotropic subundle which is integrable with respect to the $H$-twisted Courant bracket.  We discuss examples of these in \S \ref{twisted} , some of which are generalizations of $H$-twisted Poisson structures \cite{swe}.

\section{Linear Dirac Maps}\label{13Dec1}
\begin{Def} A \emph{linear Dirac structure} is a maximal isotropic subspace $D$ of $V \oplus V^*$ for some vector space $V$.  
\end{Def}
\noi In order to define morphisms between manifolds with Dirac structures,
we first develop this notion for linear maps between vector spaces equipped with linear
Dirac structures.  \\

\noindent For a vector space V, we will denote by
$p_V : V\oplus V^* \longrightarrow V $
and  $p_{V^*} : V\oplus V^* \longrightarrow V^* $
the
natural projection maps.  There are natural notions of pullbacks and pushforwards of Dirac structures, first introduced by Weinstein. For a linear map $f : V \lra W$ 
the pullback of 
a linear Dirac structure $L_W \subset  W\oplus W^*$ 
is defined as 
$f^\star L_W := \{ X + f^*\xi \in V \oplus V^* \st fX + \xi \in L \}$, 
and the pushforward of 
a linear Dirac structure $L_V$ on $V$
is $f_\star L_V := (f^*)^\star L_V$.  These are again linear Dirac structures. 
\\

\begin{Def}\label{a1} Let $L_{1},L_{2}$ be linear Dirac structures
on vector spaces ${V_1},{V_2}$.
A linear map $f: {V_1} \longrightarrow {V_2} $ is said to be a
  \emph{linear Dirac map} if:\\
\begin{equation}
  f(L_{{1}} \cap {V_1} ) \subset L_{{2}} \cap {V_2}, \; and
\tag{M1}
\end{equation}
\begin{equation}
 p_{{V_1} ^*} \inv (f^* (p_{{V_2} ^*} L_{{2}})) \cap L_{{1}} \subset f^\star L_{{2}}.
\tag{M2}
\end{equation}

\noindent Such a map will also be written as $f: ({V_1},L_{V_1} ) \longrightarrow
({V_2},L_{V_2}) $.  Notice that a linear Dirac structure on ${V_1}$ is also a linear Dirac structure 
on ${V_1}^*$, however Dirac maps are not invariant under duality.
We will say that $f$ is a \emph{linear dual-Dirac map} if
$f^*:\ ({V_2}^*,L_{V_2})
\lra ({V_1}^*,L_{V_1})$
is a linear Dirac map.

\end{Def}

\begin{lem} For a map $F: (V_1,L_1) \longrightarrow (V_2,L_2)$,
  condition (M2) of Definition~\ref{a1} is equivalent to either of
  the two following conditions:
\begin{equation}
  Y+\xi \in L_2 \; \; \mbox{and} \; \;  X + f^*\xi \in L_1 \;
  \Longrightarrow \; fX + \xi \in  L_2 ,\tag{M2$^{\prime}$}
\end{equation}
\begin{equation}
f^\star L_2 \subset L_1 + f\inv (L_2 \cap V_2) . \tag{M2$^{\prime\prime}$}
\end{equation}
\end{lem}

\begin{pf} (M$2^\prime$) is simply a rewriting of (M2) in terms of elements.  Now, (M$2^{\prime \prime}$) is equivalent to the following statements:
\[(L_1 + f\inv (L_2 \cap V_2))^\perp \subset (f^\star L_2 )^\perp = f^\star L_2,
\]
\[
L_1 ^\perp \cap (f\inv (L_2 \cap V_2))^\perp \subset f^\star L_2 ,
\]
\[
L_1 \cap (f\inv (L_2 \cap V_2))^\perp \subset f^\star L_2 .
\]
\noi So to show that (M2) and (M$2^{\prime \prime}$) are equivalent, it suffices to show that
 $(f\inv (L_2 \cap V_2))^\perp = V_1 \oplus (f^* (p_{{V_2} ^*} L_{{2}}))$
 or equivalently that $Ann(f\inv U) = f^* (p_{V_2 ^*} L_2)$.
  It is a general fact that for a map $V_1 \stackrel{f}{\lra} V_2 $ 
of vector spaces, and a subspace $U \subset V_2$,
  $Ann(f\inv U)) = f^* Ann(U)$.  
Thus, to show that $Ann(f\inv L_2 \cap V_2) = f^* (p_{{V_2} ^*} L_2)$,
 we need only show that $Ann(L_2 \cap V_2 ) = p_{V_2 ^*} L_{2}$. 
We want, therefore, that $(L_2 \cap V_2 )^\perp \subset V \oplus p_{V_2 ^*} L_2$,
 but we can rewrite this as $L_2 + V_2 = L_2^\perp + V_2 ^\perp \subset V_2 \oplus p_{V_2 ^*} L_2$,
 which is apparent.  $\square$
\end{pf}
\noi 
Recall that for any
 linear Dirac structure
$L_V$  there is a unique data of a subspace $E \subset
V$ and $\vep \in \wedge ^2 E^* $ such that $L_V$ is of the form
$
L(E, \vep ) := \{ x + \xi \in E \oplus V^* \st 
\xi_{|_E} =
\vep(x, -) \}$\ (~\cite{gua}).  
Dually, there is
a subspace $U \subset
V^*$ and  $\pi \in \wedge ^2 U^* $ such that $L_V = L(\pi , U ):=
\{x + \xi \in V \oplus U | x_{|_U} = \pi(\xi, -) \} $.
We denote the inclusion maps by 
 $j_E : E \hookrightarrow V $ and 
$i_U : U
\hookrightarrow V^*$.  \\


\noindent For
Dirac structures $L_k = L(\pi _k , U_k ) $ on $V_k$ ($k=1,2$),
consider maps $f:V_1\lra V_2$  such that $f^* (U_2 ) \subset U_1 $.  
Denote the restriction of $f^*$ to $U_2$ by  $\phi = \phi _f : U_2 \lra U_1$.  \\


\noindent 
The next two propositions~\ref{a2} and \ref{26nov1} state conditions (M1) and (M2) of Definition~\ref{a1} in terms of presentations of linear Dirac structures as $L(\pi , U)$ or $L(E,\vep)$.

\begin{prop}\label{a2} 
Let 
$L_k = L(\pi _k ,
  U_k )$
be  Dirac structures on $V_k,\ k=1,2$.  
A map  $f: V_1 \lra V_2
  $ is a  Dirac map if and only if\\
(D1) $f^* (U_2 ) \subset U_1 $
and the corresponding map $\phi=\phi_f$ satisfies\\
(D2) $\phi^* \circ \pionep \circ \phi= \pitwop $\ \
  (or equivalently, $\phi^* \pi_1 = \pi _2 $).
\end{prop} 

\noi (D1) and (D2) can be expressed by the requirement that the following diagram commutes.
\[
\begin{CD}
U_1 @<\phi _f << U_2 \\
@VV (\pi _1)_\sharp V         @VV (\pi _2 )_\sharp V\\
U_1 ^* @> \phi _f ^* >> U_2 ^*
\end{CD}
\]
\begin{pf}
Conditions (D1) and (D2) of this proposition correspond to the conditions 
(M1) and  (M2) of Definition~\ref{a1}.
\begin{enumerate}
\item We show that $f^* (U_2 ) \subset U_1 $ if and only if $f(V_1
  \cap L_1 ) \subset V_2 \cap L_2 $.  But since $f^* (U_2) \subset U_1
  $ if and only if $f(Ann(U_1)) \subset Ann(U_2) $, it suffices to
  show that $Ann(U_i ) = L_i \cap V_i $.\\

\noindent Clearly $Ann(U_i) \subset L_i \cap V_i $ by definition of
$L(\pi_i , U_i ) $.  On the other hand, $L_i \cap V_i \subset L_i =
L_i ^\perp $. Hence $\langle L_i \cap V_i , L_i \rangle = 0 $ and so $\langle L_i \cap V_i
, U_i \rangle = 0 $.  This implies that $L_i \cap V_i \subset Ann(U_i) $ and
therefore $L_i \cap V_i = Ann(U_i) $.  

\item 
First suppose that  $\phi _f ^* \circ \pionep \circ \phi _f = \pitwop
  $.  By definition, $Y+\xi \in  L_2$ if and only if $\pitwop \xi =
  i_2 ^* Y $.  Let $Y + \xi \in L_2 $ and $X+f^*\xi \in L_1$.  Then
  $\i_1 ^* X = \pionep (f^* \xi )= \pionep (\phi( \xi) ) $.  We have 
\[ \pitwop \xi = \phi ^* \circ \pionep \circ \phi (\xi) =\phi ^* \circ i_1 ^* X =
i_2 ^* \circ f (X). \]
But $i_2 ^* \circ f (X) = \pitwop (\xi )$ means exactly that $fX + \xi
\in L_2$.  Therefore condition (M2$^\prime$) is satisfied.  \\ 

\noindent 
Now suppose that  conditions (M2$^\prime$) and (M1) are satisfied.  Let $\xi \in
U_2$.  Since $f^* \xi \in U_1 = p_{V_1 ^*}L_1$, there is some $X$ such that $X + f^*
\xi \in L_1 $.  This means that $i_1 ^* X = \pionep (\phi \xi ) $.
But by condition (M2$^\prime$), $fX + \xi \in L_2 $, so  $\i_2 ^* \circ f X =
\pitwop \xi $.  However,  $i_2 ^* \circ f X = \phi ^* \circ i_1 ^* X =
\phi^* \circ \pionep \circ \phi \xi $.   Therefore,
$\phi^* \circ \pionep \circ \phi \xi = \pitwop \xi $.  But
this is true for arbitrary $\xi \in U_2$, which means $\phi^* \circ \pionep \circ \phi  =
\pitwop$.$\square$
\end{enumerate}
\end{pf}

\begin{prop}\label{26nov1} 
A linear map $f : (V_1 , L_1 = L(E_1 ,\vep _1 ) ) \lra
  (V_2 , L_2 = L(E_2 ,\vep _2 ) )$ is a linear Dirac map if and only
  if:\\
1. $f(Ker((\vep_1)_\sharp) ) \subset Ker((\vep_ 2 )_\sharp ) $, and \\
2. for $X_k \in E_k$, and $\xi_2 \in V_2 ^* $ such that 
$j_2 ^* 
\xi_2 = (\vep_ 2 )_\sharp (X_2) $ and 
$(\vep_1 )_\sharp (X_1) = j_1 ^* f^* \xi$, one has
$$	
f(X_1) \in E_2\ \ \text{and}\ \  f(X_1)-X_2\in\ Ker((\vep_ 2)_\sharp).
$$	
\end{prop} 

\begin{pf} Upon observing that $Ker((\vep_ i)_\sharp )=L_{V_i} \cap V_i$,
  (1) and (M1) are apparently equivalent. Condition (2) is a direct restatement of (M$2^\prime$) in terms of $E$ and $\vep$. $\square$
\end{pf}

\begin{prop}\label{a3} 
Pairs $(V,L) $ of vector spaces
with linear Dirac structures form a category in two ways
by taking morphisms to be either
the
  linear Dirac maps (this category will
  be denoted $\mathcal {LD} $)
or
the linear dual-Dirac maps
(category $\mathcal {LD}^* $).
The assignment $(V,L) \mapsto (V^*,L)$ gives an equivalence 
of $\mathcal {LD}^{opp}$ with $\mathcal {LD}^*$ 
\end{prop} 

\begin{pf} We must show that the composition of two linear Dirac
  maps $(V_1 , L_1 ) \stackrel{f_1}{\lra}
  (V_2 , L_2 ) \stackrel{f_2}{\lra} (V_3 , L_3 ) $ 
is a linear Dirac map.    Let
  $L_j = L(\pi _j , U_j )$, $i_j : U_j \hookrightarrow V_j $, and
  $\phi _j : U_{j+1} \lra U_j $.  Here we use the criteria of
  Proposition~\ref{a2}.\\

\noindent Because $f_1$ and $f_2$ satisfy (D1)
so does the composition:
\[
(f_2 \circ f_1 )^* (U_3 ) = f_1 ^* \circ f_2 ^* (U_3 )
\subset f_1 ^* (U_2 ) \subset U_1
\]
Again,
because $f_1$ and $f_2$ also satisfy (D2), i.e.,
$
\phi _1 ^*
  \circ \pionep  \circ \phi _1 =  \pitwop$ and $\phi _2 ^* \circ \pitwop \circ \phi _2 = (\pi_ 3)_\sharp$,
we get
\[
(\pi_ 3)_\sharp = \phi _2 ^* \circ \phi _1 ^*
  \circ \pionep \circ \phi _1 \circ \phi _2  = (\phi _1 \circ
  \phi _2 ) \circ \pionep \circ (\phi _1 \circ \phi _2 ) 
.
\]
\noindent Finally, associativity
and the existence of identity morphisms is obvious
since linear Dirac maps are functions.$\square$  
\end{pf}

\section{Dirac and Dual-Dirac Maps}\label{13Dec2}

\noi Now linear Dirac and linear dual-Dirac maps are used to define maps between manifolds with Dirac structures.

\subsection{Dirac Maps}\label{10feb2}

\begin{Def}\label{1feb1} Let $f : M\lra N $ be a $\mathcal C ^\infty $ map of
  manifolds, and let $L_M$ be an (almost) Dirac structure on M and $L_N $ an (almost)
  Dirac structure on N.  The map $f$ is said to be
  \emph{Dirac} or a \emph{Dirac map} if at each point $p \in
  M $, $df_p : (T_p M _\cc, (L_M ) _p ) \lra (T_{f(p)} N _\cc , (L_N )_{f(p)} )
  $ is a linear Dirac map. If $L_M$ and $L_N$ are both generalized complex structures, we also say that $f$ is a \emph{generalized complex map}.  A Dirac map $f$ between manifolds with Dirac
structures will also be denoted by
 $f : (M, L_M ) \lra (N, L_N ) $.
\end{Def}


\begin{prop}\label{24nov1} We can define a category $\mathcal D$ of Dirac manifolds by taking its objects to be all pairs $(M,L)$, where M is a $\mathcal C ^\infty$ manifold and L is a Dirac structure; the morphisms of $\mathcal D$ are the Dirac
maps.
\end{prop}    

\begin{pf} This follows from
  Proposition~\ref{a3}.  
For instance if
$(M_1,L_1)\stackrel{f}{\longrightarrow} (M_2,L_2)
  \stackrel{g}{\longrightarrow} (M_3, L_3) $ are Dirac maps, then for any 
$p \in M_1$ the differential
$d(g\circ
  f)_p = dg_{f(p)} \circ df_p$ is 
composition of linear Dirac maps, so it is itself a linear Dirac map.$\square$
\end{pf}

\begin{prop} Category $\mathcal D$
contains complex manifolds and Poisson manifolds as
full subcategories.  The category $\DD$ also contains the category of manifolds as a full subcategory.  
For two symplectic manifolds
$(M,\omega_M)$
and $(N,\omega_N)$
a map $f:M\to N$ which is a local isomorphism (i.e., 
 $df_p $ is an isomorphism for all $p \in M$)
is a symplectomorphism  if and only if
  f is Dirac.  
\end{prop}

\begin{pf} Consider Dirac manifolds $(M_1,L_1)$, $(M_2,L_2)$ 
and a map 
  $f :M_1 \longrightarrow M_2 $.
Since the property of being  holomorphic, Poisson, or Dirac
map  is determined pointwise, 
we fix a point $p \in M_1 $ and $q= f(p)\in M_2$.  We use the notation of Proposition~\ref{a2}.\\
 
\noi If $L_1$ and $L_2$ are complex strucutres, we need to check that
f is holomorphic if and only if f is
Dirac.
Let $E_1$ and $E_2$ be the holomorphic tangent bundle for
  complex structures so that $L_k = E_k\oplus Ann(E_k) =
  L(0,Ann(E_k))$. Because
  $\pi_i = 0 $ for $i=1,2$, condition (M2) is trivially true.  The
  map $f$ is
  holomorphic if and only if $df_p (E_1)_p \subset (E_2)_q $.  Also, $f$ is
  Dirac if and only if $df^*(Ann(E_1)) \subset Ann(E_1) $.  It is
  clear now that $f$ is Dirac exactly when $f$ is holomorphic.\\

\noi We check that if $L_1$ and $L_2$ are Poisson structures, then $f$ is a Poisson
  map if and only if f is Dirac.
Now let $\pi_1 \in \wedge ^2 T_pM_1$, $\pi_2 \in \wedge^2 T_q M_2$ be two
  Poisson bivectors with  Dirac structures 
$L_k = L(\pi _k ,T^*M_k ) $.  
Since $U_1 = T^* _p M_1$, (M1) is trivially satisfied.  
The map $f$ is Poisson if and only if 
$df _p (\pi_1 )_p
  = (\pi _2 )_q $. 
However, one may easily verify that $df _p (\pi _1 )_p = (\pi _2)_q $ if and only if $((\pi _2)_q)_\sharp = df_p \circ ((\pi _1)_p)_\sharp  \circ df_p ^*$.  Therefore, by Proposition~\ref{a2}, $f$ is Poisson if and only if $f$ is Dirac.\\


\noi Sending a manifold $M$ to $(M, TM) \in \DD$ gives a full embedding of the category of manifolds into $\DD$ since any $(X, L) \lra (M, TM)$ is Dirac for any map of manifolds $f: X \lra M$ and any Dirac structure $L$ on $X$.  \\

\noi Symplectic structures $\omega _i$ determine Poisson structures 
$\pi _i$ by $(\pi _i) _\sharp  = (\omega _i ) _\sharp \inv $ for $i =1,2$.
The claim follows from the observation 
that, since $df$
 is everywhere an isomorphism, $f$ is 
symplectic precisely if
it is Poisson, i.e., if it is Dirac. $\square$   

\end{pf} 


\begin{rem} We note that
for a Dirac map
$f : (M , L_M ) \lra (N, L_N) $, even if it is an immersion 
or a submersion,
it is not true that either of the Dirac structures
$L_M$ or $L_N$
determines the other.  For example, if $L _N = TN$ and $L_M = ( 0, U_1)$ for any subbundle $U_1$ of $TM$, then any $f$ is a Dirac map.  Similarly, any map $f$ is Dirac if $L_M = T^*M$ and $L_N = L(0,U_2)$ for any subbundle $U_2$ of $T^*N$. 
\end{rem}

\subsection{Dual-Dirac Maps}\label{13Dec3}

\begin{Def}\label{4feb1} A map $f: M \lra N$ of manifolds with (almost)
  Dirac structures $ L_M, L_N$ is said to be a \emph{dual-Dirac map} if
  $df_p$ is a linear dual-Dirac map for each $p \in M$.   
\end{Def} 

\noi The dual statement of Proposition~\ref{a2} is as follows.

\begin{prop}\label{24mar} Let $V_1$, $V_2$ be vector spaces and $L_1 = L(E_1 ,
  \vep _1 ) \subset V_1 \oplus V_1 ^* $, $L_2 = L(E _2 , \vep _2 ) \subset
  V_2 \oplus V_2 ^* $ be linear Dirac structures.  If $f: V_1 \lra V_2
  $ is any linear map, then f is a linear dual-Dirac map if and only if\\
$(D1^* )$  $f (E_1 ) \subset E_2 $ and \\
$(D2^* )$  $f  \circ (\vep _2)_\sharp \circ f^* = (\vep_1)_\sharp $ (or equivalently, $f^* \vep _2 = \vep _1 $).

\end{prop} 
\noi These two conditions are equivalent to the existence and commutivity of the diagram 
\[
\begin{CD}
E_1 @>> f> E_2 \\
@VV (\vep _1)_\sharp V   @VV (\vep _2)_\sharp V \\
E_1 ^* @<< f^* < E_2 ^*
\end{CD}
\]

\begin{prop} Dirac manifolds with dual-Dirac maps form
  a category $\mathcal{D}^*$ of dual-Dirac manifolds.
\end{prop}
\begin{pf}
If $U \stackrel{f}{\lra} V \stackrel {g}{\lra} W$ are two linear dual-Dirac maps, then $(g \circ f)^* = f^* \circ g^*$ is the composition of two linear dual-Dirac maps, so $g \circ f$ is linear dual-Dirac.  Now the conclusion follows in the same way as for Proposition~\ref{24nov1}.$\square$
\end{pf}

 \begin{prop} The category $\DD ^*$ of dual-Dirac manifolds contains symplectic and complex manifolds as full subcategories.  
\end{prop}

\begin{pf} If $(M_1,L_1) , (M_2,D_2) \in \DD ^*$ and $L_1$ and $L_2$ are complex structures, then $E_1$ and $E_2 $ are the holomorphic tangent bundles, and $\vep _1 = \vep _2 =0$.  In light of Proposition~\ref{24mar}, a map $f : M_1 \lra M_2$ is holomorphic if and only if it is dual-Dirac.  If $L_j = L(TM_j , i\vep _j)$ ($j = 1,2$) are symplectic structures, then Proposition~\ref{24mar} states that $f :M_1 \lra M_2$ is dual-Dirac if and only if $f$ is a symplectomorphism.  $\square$
\end{pf}

\begin{rem} 
We say that a Dirac manifold $(M,L_1)$ is a \emph{(dual)Dirac submanifold} of $(N,L_2)$ if $M $ is a submanifold of $N$ and the inclusion map $i: (M, L_1) \hookrightarrow (N , L_2)$ is a (dual)Dirac map.  
If both Dirac structures are generalized complex structures
this gives two notions of generalized complex submanifolds.  
Gualtieri~\cite{gua} offers a very different definition of generalized 
complex submanifold which is related to the notion of \emph{branes}.  The terminology ``Dirac submanifolds'' follows the standard use of the notion of a subobject in a category (that the inclusion is a morphism in the category).  Gualtieri's ``generalized complex submanifolds'' are submanifolds with extra structure which formalize mathematically certain subclasses of branes, objects from physics which are viewed as boundary conditions for Quantum Field Theories.  These are objects of a different nature than the ambient manifold and are in particular not categorical subobjects.  The inclusion maps of Ben-Bassat and Boyarchenko \cite{bbb} can also be used to define a type of generalized complex submanifolds.
\end{rem}



\subsection{B-transforms} 
\noi For a vector space $V$ and $B \in \wedge ^2 V^*$, we extend the composition $ V \stackrel{B_\sharp}{\lra} V^* \hookrightarrow V \oplus V^*$ by $0$ on $V$ to give a map $\tilde B _\sharp : V \oplus V^* \lra V \oplus V^*$ and a map $e ^{\tilde B _\sharp}$, which by abuse of notation, we denote by $e ^B$.  This is the linear version of a B-transform, but we also refer to it as a B-transform.  For a linear Dirac structure $L \subset V \oplus V^*$, we say that $e^BL$ is the \emph{B-transform of L} by $B$.  If $L = L(E  , \vep)$, then $e ^B L = L(E, \vep + B_{|E\times E})$.
 
\begin{prop} Dual-Dirac maps are stable under B-transforms in the
  sense that if $B \in \Omega ^2 (N)$, then  $f : (M, L_M) \lra (N ,
  L_N)$ being dual-Dirac implies $f : (M, e^{f^*B}L_M) \lra (N ,
  e^BL_N)$ is dual-Dirac. 
\end{prop}

\begin{pf} For dual-Dirac maps, it is enough to prove this statement
  pointwise for the derivative $df$ of $f$, where one may represent
  $L_M$ as $L_M = L(E_1 , \vep _1)$ and  $L_N$ as $L_N = L(E_2 , \vep
  _2)$.  Clearly $f^* \vep _2 = \vep _1$ implies $(f_{|E_1})^*(\vep _2 + B_{|E_2}
  ) = \vep _1 + f^*B _{|E_1}$.$\square$  
\end{pf}  

\begin{rem} Crainic's notion ~\cite{cra} of generalized complex maps is
  also stable under $B$-transforms.  However, Dirac maps are not stable under B-transforms, as the following example shows.  Let $\io : M \hookrightarrow M\times N$ be any inclusion $x \mapsto (x,c)$, let $L_1 = TM$, and let $L_2 = T(M\times N)$.  For any closed non-zero $B \in T^*M \wedge T^*N \subset \wedge ^2 T^* (M\times N)$, $\io ^* B = 0$, and $(M , e^0L_1) \lra (M\times N , L( T(M\times N) , B))$ is not Dirac.  
\end{rem}

\subsection{Comparison with Other Categories of Dirac Manifolds}\label{comp}
\noi In addition to the Dirac maps and dual-Dirac maps described in this paper, there are several other distinct concepts of Dirac maps presented in the existing literature, some of them only defined for a subclass of Dirac manifolds.  \\

\noi Crainic~\cite{cra} defines a category of generalized complex manifolds which also contains complex manifolds as a full subcategory but which is only defined for generalized complex structures, not for Dirac structures in general. He defines ``generalized holomorphic maps," which visibly form a category.  The requirement to be a \emph{Crainic-Dirac} map is very strong.  For instance, if $L_M$ and $L_N$ are symplectic structures, then $f : (M, L_M) \lra (N, L_N)$ is Crainic-Dirac if and only if it is both symplectic and Poisson.   Recall from \S \ref{13Dec1} that there are pushforwards and pullbacks of linear Dirac structures.  Both pushforwards and pullbacks are useful operations in their own right, but they also define maps between manifolds with Dirac structures \cite{bur}.  One asks of $f: (M,L_M) \lra (N, L_N)$ that $L_N = f_\star L_M$ or $L_M = f^\star L_N$.  However, when considered as maps of manifolds with Dirac structures, they do not generalize holomorphic maps.  Pullbacks and pushforwards only generalize holmorphic maps when $f$ is an immersion or submersion respectively.  \\

\noi Alekseev, Bursztyn, and Meinrenken \cite{abm} have introduced a notion of Dirac maps, which consists of a pair $(f,B)$ of a map of manifolds $f : M \lra N$ and a 2-form $B \in \Om ^2 (M)$ such that $(f,B): (M,D) \lra (N, L)$ satisfies $f_\star (e^B D) = L$, where $f_\star L$ is the pushforward of $L$.  This fits naturally into the picture when one considers generalized complex structures from the point of view of pure spinors.  It is also useful for dealing with Courant algebroids other than $\VV _M$.  Our notions (defined below) of maps of Courant algebroids and ``Courant-Dirac maps'' use the same data $(f,B)$ of a map and a 2-form as in \cite{abm}. \\

\noi Additionally, Ben-Bassat and Boyarchenko \cite{bbb} define Dirac inclusions and quotients which combine to give a notion of Dirac maps distinct from ours.  Most recently, Ornea and Pantilie have defined ``generalized holomorphic maps''~\cite{op}, which they have used to study generalized K\"ahler structures and which are compatible with Poisson and complex maps as well as being stable under B-transforms.

%
%
%


\section{Cateogries of Courant Algebroids and Dirac Manifolds}\label{10feb1}

\noi Here we extend the definitions of Dirac and dual-Dirac maps to include Dirac structures contained in arbitrary exact Courant algebroids.

\subsection{Exact Courant Algebroids}
 
 \noi A Courant algebroid on $M$ (defined in \cite{lwx}) is a quadruple $(\EE, \pi : \EE \lra TM , \langle \, , \, \rangle , [ \, , \, ])$, where $\langle \, , \, \rangle$ 
is bilinear form and 
$[ \, , \, ]$ is a bracket operation on sections of $\EE$. If we define $D : \CC ^\infty (M) \lra \Ga (M, \EE)$ by $\langle Df , A \rangle = \frac {1}{2} \pi (A) f$.  We say that $(\EE, \pi : \EE \lra TM , \langle \, , \, \rangle , [ \, , \, ])$ is a \emph{Courant algebroid } if for all sections $A,B,C$ of $\EE$ and functions $f, g$ on $M$, the following conditions hold: \\
1. $\pi ([A,B])= [\pi A , \pi B]$\\
2. $Jac (A,B,C) = D (Nij(A,B,C))$\\
3. $[A,fB] = f[A,B] + (2\langle Df ,A \rangle )B - \langle A , B \rangle Df $\\
4. $\pi \circ D = 0$ \\
5. $\pi (A) \langle B , C \rangle = \langle [A,B] + D\langle A , B \rangle , C \rangle + \langle B , [A,C] + D \langle A , C \rangle \rangle $.\\
\noi $Jac$ denotes the Jacobiator of the Courant bracket, and $Nij (A,B,C) := \tfrac{1}{3}(\langle [A,B],C \rangle + \langle [B,C],A \rangle + \langle [C,A], B \rangle ) $.\\

\noi An \emph{exact Courant algebroid}
 is a Courant algebroid $\EE \lra M$
 for which $0 \lra T^* M \stackrel{\pi ^* }{\lra} \EE \stackrel{\pi}{\lra} TM \lra 0$ 
is exact.  A map of Courant algebroids on $M$ is simply a map of vector bundles which preserves all of the defining structures of the Courant algebroids.  Such a map is necessarily an isomorphism if the Courant algebroids are exact.  For a closed 3-form $H$, there is an exact Courant algebroid 
$\VV _{M,H}$, where 
$(\VV _{M,H} , \pi , \langle \, , \, \rangle ) = (\VV _{M} , \pi , \langle \, , \, \rangle )$,
 but $\VV _{M,H} $
 has bracket given by $[X + \xi , Y + \eta ]_H = [X+\xi , Y + \eta ] + \io _Y \io _X H$ 
\cite{swe}.  \v Severa's classification (as described in \cite{bcg}) is that any exact Courant algebroid is isomorphic to some $\VV _{M,H}$. 
 Further, it is known that any isomorphism $\VV _{M,H} \lra \VV _{M,K}$ 
is a B-transformation $e^B$, where $H= K + dB$~\cite{gua}. \\

\noi We now define the pullback of an exact Courant algebroid.  For a submanifold $\io :  S \hookrightarrow M$ and an exact Courant algebroid $\EE$ over $M$, $\EE$ can be restricted to an exact Courant algebroid $\EE _{|S} = \pi \inv (TS) / Ann(TS)$ on $S$ \cite{bcg}, and when $\EE = \VV_{M,H}$, there is a canonical identification $(\VV _{M,H})_{|S} = \VV _{S, \io ^* H}$.  For a map of manifolds $f: M \lra N$ and a Courant algebroid $\EE$ over $N$, we define the pullback $f^\star \EE$ as follows.  The restriction of $\VV _M \oplus \EE$ to the graph $\Ga _f$ of $f$ is an exact Courant algebroid on $\Ga _f$.  The diffeomorphism $F = 1 \times f$ between $M$ and $\Ga _f$ now gives an exact Courant algebroid on $M$, which will be denoted by $f^\star\EE$.  Clearly, $f^\star \EE = (\VV _M \oplus _{T_M} f^* \EE )/F$, where $f^* \EE$ is the pullback of vector bundles and $F$ is the subbundle of $\VV _M \oplus _{TM} f^*\EE$ which is the graph of $-f^* : f^* T^*N \lra T^*M$.  Concretely, 
at a point $x \in M$, the fiber of $f^\star \EE$ is 

\[ (f^\star \EE )_x = \{ (A,B) \in (\VV _M )_x \oplus \EE _{f(x)} \st \pi B  = d_x f \pi A \} / \{ (-f^*\be , \be) \st \be \in T^* _x N\}.
\] 
\noi When $\EE = \VV _N$, there is an isomorphism $\VV _M \lra f^\star \VV _N$ given by $X + \al \mapsto [X + \al , df X ]$, and more generally when $\EE = \VV _{N,H}$, $f^\star \EE $ is identified with $\VV _{M, f^*H} $ in the same way because $\EE _{|\Ga _f} = \VV _{M\times N , 0+H}$.\\

\begin{lem} Let $\EE$ be an exact Courant algebroid on $M$.  When $\io : S \hookrightarrow M$ is an embedding, $\io ^\star \EE \simeq \EE _{|S}$.  For any manifold $N$, let $p : M \times N \lra M$ be the first projection map.  Then $p ^\star \EE \simeq \EE \oplus \VV _N$.
\end{lem}
\begin{pf} 
The pullback $\io ^\star \EE$ is the exact Courant algebroid $(\VV _S \oplus \EE)_{|\Delta (S)} \lra \Delta S \tilde \lra S$, where $\Delta : S \lra S \times S$ is the diagonal map.  There is an isomorphism $\EE _{|S} \lra (\VV _S \oplus \EE )_{\Delta S}$ given by $[v] \mapsto [\pi v , v]$.  For the projection map $p$, $p ^\star \EE = (\VV _M \oplus \VV _N \oplus \EE )_{|N \times \Delta M} = \VV _N \oplus (\VV _M \oplus \EE )_{|\Delta M} = \VV _N \oplus id_M ^\star \EE $.  But $\VV _N \oplus id_M ^\star \EE = \VV _N \oplus \EE$ by Lemma~\ref{11feb1} below. $\square$
\end{pf}

\subsection{Category of Spaces with Courant Algebroids}
\noi A \emph{space with a Courant algebroid} is simply a Courant algebroid $\EE \lra M$ over a manifold $M$, also denoted $(M, \EE)$.  A morphism  $(M, \EE ) \lra (N , \EE ^\prime)$ between two spaces with Courant algebroids is a pair $(f, \phi)$ of a map $f : M\lra N$ together with an isomorphism of Courant algebroids $\phi : \EE \lra f^\star \EE ^\prime$. 

\begin{lem}\label{11feb1} The collection of spaces with Courant algebroids and the morphisms between them form a category $\CC$.
\end{lem}

\begin{pf} The only thing to check is that morphisms compose.  Let $(M_1 , \EE _1 ) \stackrel {(f,\phi)}{\lra} (M_2 , \EE _2 ) \stackrel {(g,\psi)}{\lra} (M_3 , \EE _3)$.  
We first observe that for any isomorphism $\al : \EE \lra \EE ^\prime$
 of exact Courant algebroids on manifold $M$ and any submanifold $\io : S \hookrightarrow M$, $\al $ restricts to an isomorphism $\al _{|S} : \EE _{|S} \lra \EE ^\prime _{|S}$. 
 This implies that for any map $f: N \lra M$, there is a pullback
 $f^\star \al : f^\star \EE \lra f^\star \EE ^\prime$.
  Since there are isomorphisms $\EE _1 \stackrel {\phi}{\lra} f^\star \EE _2$ 
and $\EE _2 \stackrel {\psi}{\lra} g^\star \EE _3$,
 there is also an isomorphism $ f^\star \EE _2  \stackrel {f^\star \psi} {\lra} f^\star g^\star \EE _3$.
  Therefore, there is an isomorphism $f^\star \psi \circ \phi : \EE _1 \lra f^\star g^\star \EE _3$.
  To complete the proof, we must show that $f^\star g^\star \EE _3 = (g\circ f )^\star \EE _3$. \\

\noi Let $\al : M_1 \tilde \lra \Ga _f$ and $\be : M_2 \tilde \lra \Ga _g$.  There is a commutative diagram 
\[ 
\begin{CD} 
M_1  @> \al >>  \Ga _f @>i>> M_1 \times M_2  @. \\
@|  @VV\ga V  @VV \delta V @. \\
M_1 @> \epsilon >>  \Ga _{\be \circ f} @>j>> M_1 \times \Ga _g @>k>> M_1 \times M_2 \times M_3    
\end{CD}
\]

\noi where $\al$, $\be$, $\delta$, $\ga$, $\epsilon$ are isomorphisms, and $i$, $j$, and $k$ are inclusions.  Starting with $\VV _{M_1} \oplus \VV _{M_2} \oplus \EE _3$ on $M_1 \times M_2 \times M_3$, we restrict the Courant algebroid for each inclusion and make identifications (pullback of vector bundles) for each isomorphism.  Along the top row, we end up with $f^\star g^\star \EE_3$ on $M_1$, and the result on the second row will therefore be the same.  \\

\noi Next we observe that the inclusion $\Ga _{\be \circ f} \hookrightarrow M_1 \times \Ga _g \hookrightarrow M_1 \times M_2 \times M_3$ is the same as the inclusion $\Ga _{\be \circ f} \hookrightarrow \Ga _f \times M_3 \hookrightarrow M_1 \times M_2 \times M_3$.  In general, for $Q \subset S \subset M$ and Courant algebroid $\EE$ on $M$, $\EE _{|Q} = (\EE _{|S})_{|Q}$.  Hence, restricting $\VV _{M_1} \oplus \VV _{M_2} \oplus \EE _3$ to $\Ga _{\be \circ f}$ is the same as restricting $\VV _{\Ga _f } \oplus \EE _3$ on $\Ga _f \times M_3$ to $\Ga _{\be \circ f}$.  There is now a commutative diagram 
\[ 
\begin{CD} 
M_1  @> \al >>  \Ga _{\be \circ f} @>l>> \Ga _f \times M_3   \\
@|  @VV\io V  @VV \al \inv \times 1 V  \\
M_1 @> \kappa >>  \Ga _{g \circ f} @>m>> M_1 \times M_3     
\end{CD}
\]

\noi where $\io$ and $\kappa$ are isomorphism and $l$ and $m$ are inclusions.  Again restricting for each inclusion  and making identifications for each isomorphism, the bottom row will yield $(g\circ f)^\star  \EE _3$ on $M_1$, whereas the previous diagram shows that the top row yields $f^\star g^\star \EE _3$.  $\square$

\end{pf}

\subsubsection{Categories of Courant-Dirac manifolds}\label{twisted} 

\noi First consider Dirac structures $D_1 \subset \VV _{M_1 , H_1}$, $D_2 \subset \VV _{M_2 , H_2}$.  Since $\VV _{M_i , H_i} = \VV _{M_i}$ as vector bundles, it still makes sense to say that a map $f : M_1  \lra M_2 $ is Dirac or dual-Dirac if it satisfies the conditions of Definitions~\ref{1feb1} and \ref{4feb1} respectively.  

\begin{Def}\label{4feb2} We define \emph{Courant dual-Dirac manifolds} to be triples $(M, \EE , D)$ of a Dirac structure $D$ in an exact Courant algebroid $\EE$ on manifold $M$.  A morphism of Courant dual-Dirac manifolds $ (f, \phi): (M_1 , \EE _1 , D_1 ) \lra (M_2 , \EE _2 , D_2)$ is a map $(f,\phi) : (M_1 , \EE _1 ) \lra (M_2 , \EE _2 )$ of spaces with Courant algebroids such that:\\
$(M1)^*$ : For any $[A,B] \in f^\star \EE _2$, if $[A,B] \in \phi D_1$, then there exists $\be \in T^* _{f(x)} N$ such that $B + \be \in D_2$.  \\
$(M2)^*$ : If $[X + \al , B] \in  \phi D_1$ with $B \in D_2$, then $[X, B] \in  \phi D_1$.  
\end{Def} 

\noi As an example, a Courant dual-Dirac map $(M _1, \VV _{M _1, H_1} , D_1 ) \lra (M_2 , \VV _{M _2, H_2} , D_2 )$ consists of a pair $(f, B)$ such that $H_1 + dB = f^* H_2$ and $f: (M_1, e^B D_1 ) {\lra} (M_2 , D_2)$ is dual-Dirac in the sense of Definition~\ref{4feb1}.  This follows from \v Severa's classification \cite{bcg}. 

\begin{lem} 
\begin{enumerate}
\item Courant dual-Dirac manifolds form a category, which we denote by $\CC \DD ^*$. 
\item For the full subcategory $\TT \DD ^*$ of all $(M, \VV _{M , H}, D)$, the inclusion into $\CC \DD ^*$ is an equivalence. 
\item The category $\DD ^*$ is a subcategory of $\CC \DD ^*$, obtained by allowing only $\EE$'s of the form $ \VV _M$ and only allowing maps $(f,B)$ with $B = 0$.
\end{enumerate}
\end{lem}
\begin{pf}
 Since $\DD ^*$ is a category and morphisms in $\DD ^*$ are stable under B-transforms, it follows that $\TT \DD ^*$ is a category. Morphisms in $\CC \DD ^*$ have the property that if $\al _i : \EE _i \lra \EE _i ^\prime$ are isomorphisms of Courant algebroids on $M_i$, then for $(f , \phi ) \in Hom _\CC ( (M_1, \EE _1 ), (M_2, \EE _2)$, it is the case that $(f,\phi) : (M_1 , \EE _1 , D_1 ) \lra (M_2 , \EE _2 , D_2)$ is a morphism in $\CC \DD ^*$ if and only if $(f, f^\star \al _2 \circ \phi \circ \al _1 \inv ) : (M_1 , \EE _1 ^\prime , \al _1 D_1 ) \lra (M_2 , \EE _2 ^\prime, \al _2 D_2)$ is a morphism in $\CC \DD ^*$.   Hence, since every object in $\CC \DD ^*$ is isomorphic to some object in $\TT \DD ^*$ and $\TT \DD ^*$ is a category, morphisms in $\CC \DD ^* $ compose, making it a category equivalent to $\TT \DD ^*$.  Clearly $\DD ^*$ is a subcategory. $\square$
\end{pf} 

\begin{Def} We define a collection $\CC \DD$ of \emph{Courant Dirac manifolds} where $Ob(\CC \DD) $ consists of quadruples $(M , \EE , \al , D)$ of a manifold $M$ with exact Courant algebroid $\EE$, an isomorphism $\al  : \EE \lra \VV _{M,H}$ for some $H \in \Om ^3 (M)$, and a Dirac structure $D \subset \EE$.  A morphism $ (M_1 , \EE _1 , \al _1 , D_1 ) \lra (M_2 , \EE _2 , \al _2 , D_2 )$ is a map $(f,\phi) : (M_1 , \EE _1 ) \lra (M_2 , \EE _2)$ of spaces with Courant algebroids such that  $f : (M_1 ,  \al _1 D_1 ) \lra (M_2 , \al _2 D _2)$ is a Dirac map in the sense of Definition~\ref{1feb1}.  
\end{Def}

\begin{rem} 
\begin{enumerate}
\item Since $\CC$ and $\DD$ are categories, $\CC \DD$ is a category.
\item There is no way to consistently state (M1) and (M2) (as in Definition \ref{4feb2}) for Dirac structures contained in arbitrary exact Courant algebroids.  This is possible for $(M1)^*$ and $(M2)^*$ because dual-Dirac maps are stable under B-transforms whereas Dirac maps are not.  This failure forces us to choose isomorphisms $\al : \EE _i \lra \VV _{M _i , H _i}$, or equivalently isotropic sections $TM_i \stackrel{s_i}{ \lra} \EE_i$ of $\pi _i$.  In other words, different choices of isomorphisms $\al _1$, $\al _2$ will change whether $f : (M_1 , \al _1 D_1 ) \lra (M_2 , \al _2 D _2)$ is a dual-Dirac map in the sense of Definition~\ref{1feb1}.  
\end{enumerate}

\end{rem}

\section{Dirac Groups}\label{13Dec4}

\noi Here we will consider groups in the categories $\DD$ and $\DD ^*$ but will also briefly mention groups in $\TT \DD$ and $\TT \DD ^*$.

\subsection{Almost Dirac Groups}

\noi In this subsection we describe the data and conditions involved in requiring that a group with almost Dirac structure is a Dirac group (that is, group multiplication is a Dirac map).

\begin{Def} If $(G,D)$ is a Lie group with (almost) Dirac structure $D$, we say
  that $(G,D)$ is a \emph{(almost) Dirac group} if group multiplication $\mu :
  G \times G \lra G$ is a Dirac map. If $D$ is a real Dirac structure and we say that $(G,D)$ is a \emph{real Dirac group}, whereas if $D$ is a generalized complex structure, we say that $(G,D)$ is a \emph{generalized complex group}.
\end{Def} 

\noi For a Lie group $G$ and a bi-invariant subbundle $U \subset T^* G _\cc$, a section $\be$ of $\wedge ^2 U^*$ is called \emph{multiplicative } if $d\mu _{(g,h)} (\be _g + \be _h ) = \be _{gh}$ for all $g,h \in G$, i.e. $dL_g \be _h + d R_h \be _g  = \be _{gh}$.  For the remainder of this section we will, when convenient, use left translation to identify $TG$ and $T^*G$ with $G \times \fg$ and $G \times \fg ^*$. This identifies $\Ga (G , TG) \simeq Map (G, \fg)$, and with bi-invariant $U$ as above, $\Ga (G , \wedge ^2 U^*) \simeq Map (G , \wedge ^2 U_e ^*)\simeq Map (G , \wedge ^2 \cts / Ann(U_e))$.  

\subsection{Almost Dirac Groups}

\begin{lem}\label{24nov3} Let $(G,D) $ be a group with almost Dirac structure.  Then $(G,D)$  is an almost Dirac group if and only if there is a bi-invariant subbundle $U$ of $T^*G _\cc$ and multiplicative section $\beta$ of $ \wedge ^2 U^*$ such that $D = L(\beta , U)$.  In this case $\fk = Ann(U_e )$ is a $G$-invariant ideal of $\cts$.  Almost Dirac groups are thus parameterized by pairs $(\fk, \be)$ of a G-invariant ideal $\fk \subset \cts$ and a multiplicative $\be : G \lra \wedge ^2 \cts / \fk$.  


\end{lem}

\begin{pf} Let $\mu : G \times G \lra G$ denote group multiplication
  in G.  The map $d\mu ^* _{(g,h)} : T^*_{gh} G \lra T_g ^* G \oplus
  T_h ^* G $ is given by $d\mu ^* _{(g,h)} = (dR_h)_g ^* \oplus
  (dL_g)_h ^*$.  We prove this lemma for real Dirac structures, but the proof applies to arbitrary Dirac stuctures by complexifying $\VV _G$.  \\

\noi First we show that (M1) is equivalent to the requirement that there exists a bi-invariant subbundle $U$ for which $D$ is of the form $L(\be, U)$ for some $\be \in \Ga ( G, \wedge ^2 U^*)$.  At each point $g \in G$, we can express $D_g $ as $D_g = L(\beta _g , U_g)$ for some $U _g \subset T_g ^*G$ and $\beta _g \in \wedge ^2 U_g ^*$.  This defines some $U \subset T^*G$, and we must show that $U$ is a subbundle.  The Dirac structure $\mathcal D$ on $G\times G$ is
defined by $\mathcal D _{(g,h)} = D_g \oplus D_h = L(\beta _g + \beta
_h , U_g \oplus U_h )$.  In order to satisfy (M1),
$d\mu _{(g,h)} ^* (U_{gh} ) \subset U_g \oplus U_h $, which means
$(dR_h)_g ^* U_{gh} \subset U_g $ and $(dL_g)_h ^* U_{gh}\subset U_h$.\\

\noindent This is true for all $g, h \in G$.  In particular, if $h= e$,
$(dL_g)_e ^* U_g \subset U_e$, and if $g=e$, $(dR_h)_e ^* U_h \subset
U_e$.  Hence, $ U_g \subset (dL_g)_e ^{-*} U_e $ and $U_g \subset
(dR_g)_e ^{-*} U_e$ for all $g \in U_g $.\\

\noindent On the other hand, if $h = g\inv $, we get $(dR_{g\inv} )_g ^* U_e
\subset U_g$. This implies that dim$ U_g$ $\geq$ dim$U_e $. Therefore,
$U_g = (dL_g)_e ^{-*} U_e = (dR_g)_e ^{-*} U_e$. We conclude that $U$
is a subbundle and in fact bi-invariant, which
means that $U_e$ is $G$-invariant and $\mathfrak k :=
Ann(U_e) $ is $G$-invariant.  Therefore $\mathfrak k$ is an ideal.

\noi Condition (M2) states that
$d\mu (\be _{g} + \be _{h} ) = \be _{gh } $, just as
for Poisson groups.  $\square$
\end{pf} 
 
 \begin{rem} It follows from Lemma~\ref{24nov3} that if $(G , L(\be , U))$ is a Dirac group, then $\be _e =0$ and $D_e = \mathfrak k  \oplus Ann(\mathfrak k)$, where $\fk = Ann(U_e)$.  With $TG \simeq G \times \fg$, multiplication $d\mu _{(g, h)} :
\mathfrak g \oplus \mathfrak g \lra \mathfrak g $ is given by $d\mu
_{(g, h)} (X,Y) = Ad(h\inv ) X + Y$, and $d\mu _{(g,h)} ^* \xi = Ad(h)^{-*}\xi + \xi $.  Hence $\be : G \lra \wedge ^2 \cts / \fk$ is multiplicative if and only if $\be _{gh} = \be _h + Ad(h\inv)\be _g$.  
 \end{rem}

\noi The following is a generalization of the Poisson case from~\cite{luw}, \cite{vai} of muliplicative bivectors in terms of cocycles.  For an ideal $\fk \subset \cts$, $\be : G \lra \wedge ^2 \cts / \fk$ is a \emph{cocycle} if $ad_x (d_e\be (y)) - ad_y (d_e \be (x)) - d_e \be ([x,y]) = 0 $ for all $ x, y \in \fg$.

\begin{prop}\label{5feb1}
Let $G$ be a Lie group, $\fk$ a G-invariant ideal and $\be \in \Ga (G , G \times \wedge ^2 \cts / \fk)$ such that $\be (e) = 0$.  If $\be$ is multiplicative, then $d_e \be :\fg \lra \wedge ^2 \cts / \fk$ is a $\wedge ^2 \cts / \fk$-valued cocycle.  Conversely, if $G$ is simply connected, then any cocycle $\epsilon : \cts \lra \wedge ^2 \cts / \fk$ is the differential of some multiplicative $\be : G \lra \wedge ^2 \cts /\fk$. 
 
\end{prop}
\begin{pf}  
Since $G\times \fk$ is an ideal in the sheaf of Lie algebras $\TT G$, it makes sense to take the Lie derivative $\LL_X \be$ for any left or right-invariant vector field $X$. The proof of (1) is now identical to the one in \cite{vai}.  $\square$
\end{pf}

\begin{lem} If $(G,D)$ is an almost Dirac group, then inversion $\io : G \lra G$ ($g \mapsto g\inv$) and the identity map $\{e\} \hookrightarrow G$ are Dirac maps, so $(G, D)$ is a group object in the category of almost Dirac manifolds.  
\end{lem}
\begin{pf}  
By Lemma~\ref{24nov3}, $D = L(\be , U)$, and $\be _e = 0$.  Now (D1) and (D2) are trivially satisfied for $\{e\} \hookrightarrow G$, so the identity map is Dirac.  The derivative of inversion at $g\in G $ is given by $d\io _g = -Ad(g)$.  Since $U$ is bi-invariant, by Lemma~\ref{24nov3}, (D1) is satisfied.  Since $\be$ is multiplicative, $\be _{gh} = \be _h + Ad(h\inv)\be _g$.  Letting $h = g\inv$, we get $\be _{g\inv} = -Ad(g)\be _g$, which is exactly what is needed for (D2) to be satisfied for $\io$. $\square$
\end{pf}

\begin{prop} Group multiplication $\mu : G \times G \lra G $ is Dirac if and only if $(\mu , 0)$ is an ABM-Dirac map (i.e. it is a Dirac map in the sense of \cite{abm}).
\end{prop}

\begin{pf} First we observe that a map $(f,0) : (M_1 , D_1 ) \lra (M_2 , D_2)$ is an ABM-Dirac map if and only if $X + f^* \xi \in D_1 \Longrightarrow fX + \xi \in D_2$.  Let $(G,D)$ be a group in $ABM -\DD$.  Pointwise we can represent $D$ as $D_g = L(\be _g , U_g)$.  Using left translation to get a trivialization $\VV _G \simeq G\times (\fg \oplus \fg ^*)$, we view $D _g \subset \fg \oplus \fg ^*$.  The condition for $\mu$ to be a map in $ABM-\DD$ is: for all $g,h \in G$, 
\[ X + Ad(h\inv)^* \xi \in D_g \; \;  and \; \;  Y + \xi \in D_h \Longrightarrow Ad(h\inv)X + Y + \xi \in D_{gh}. \]

\noi In particular this implies that $Ad(h\inv)^* Ann(U_g) + Ann(U_h) \subset D_{gh}$.  But since $D_{gh} \cap \fg = Ann(U_{gh})$, we have $Ad(h\inv)^* Ann(U_g) + Ann(U_h) \subset Ann(U_{gh})$.  Then of course $Ann(U_h) \subset Ann(U_{gh})$ and $U_{gh} \subset U_h$ for all $g$, $h \in G$.  Letting $h = e$ gives $U_g \subset U_e$, and letting $g = h\inv$ gives $U_e \subset U_h$ so that $U_g = U_e$ for all $g \in G$.  But $Ad(h\inv)^* Ann(U_g) + Ann(U_h) \subset Ann(U_{gh})$ also implies that $Ad(h)^*Ann(U_e) \subset Ann(U_e)$ for all $h \in G$, whence $U_e$ and $\fk = Ann(U_e)$ are $G$-invariant.  Now the condition for $\mu $ being a morphism in $ABM-\DD$ is that $(\be _{gh} )_\sharp = Ad(h\inv) (\be _g)_\sharp Ad (h\inv)^* + (\be _g)_\sharp$, i.e. $\be _{gh} = Ad(h\inv) \be _g + \be _h$, which is the same as saying that $\be$ is multiplicative.  Therefore, $(G, D)$ is a Dirac group.  Conversely, let $(G , D = L(\be , U))$ be a Dirac group.  Since $U$ is bi-invariant, $(M1)$ and $(M2 ^\prime)$ easily imply that $\mu$ is a morphism in $ABM-\DD$.  $\square$
\end{pf}

\subsection{Real Dirac Groups}
\noi Here we consider Dirac groups $(G,D)$ for which $D$ is a real Dirac structure.  We prove Proposition~\ref{4auga} which relates real Dirac group structures on $G$ to Poisson group structures on quotients of $G$.

\begin{prop}\label{4auga} Let $(G, D)$ be a real Dirac group.  If the connected subgroup $K \subset G$ with $Lie \, K = \mathfrak k$ is closed, then the data for the Dirac structure $D$ on $G$ is equivalent to the data for a Poisson group structure on $G/K$.  If $G$ is semisimple, then Dirac group structures on $G$ are in bijection with Poisson group structures on quotients of $G$ by closed, normal, connected sugroups of $G$.  
\end{prop}

\begin{pf} 
\noi Suppose that there is a closed subgroup $K$ such that $Lie \, K = \mathfrak k$. Let $\pi : G \lra G/K $ denote the quotient map. When we view $G $ locally as a product of an open set $W$ of K
and an open set $V$ of G/K, we see that $T_g K = dL_g \mathfrak k$ and so $T^*V = U$.  One may check that $\pi _\star D$ is an almost Dirac structure, and $D = TK \oplus \pi ^* (\pi _\star
D) $.  Then $\pi ^\star (\pi _\star D) = D $ implies that $\pi _\star
D$ is a Dirac structure.  Since $pr_{T^* (G/K) } \pi _\star D = T^*(G/K)$, $\pi
_\star D = L(\beta , T^* (G/K))$ for some 2-form $\beta$ on $G/K$.
But $\beta$ is a Poisson structure on $G/K$ since $\pi _\star D $
is involutive.  \\

\noi To see that $(G/K , \beta ) $ is a Poisson group, we must must
show that $\beta$ is multiplicative. Let $\mu : G\times G \lra G$ and
$\overline \mu : G/K \times G/K \lra G/K $ denote group multiplication
in $G$ and $G/K$ respectively.  We know that $\overline \mu \circ \pi
\times \pi = \pi \circ \mu$.  Recall that $\gamma$ is multiplicative, meaning $d\mu (\gamma _{g_1} + \gamma _{g_2} ) = \gamma _{g_1 g_2 } $ just as
for Poisson groups.   We also know that $\beta = (d\pi) \gamma $
because $\pi :(G,D) \lra (G/K , \pi _\star D ) $ is obviously a Dirac
map.  Therefore $\beta _{g_1K} + \beta _{g_2 K} = d(\pi \times \pi
)_{(g_1 , g_2 )} 
(\gamma _{g_1} + \gamma _{g_2} ) $, so
 \begin{eqnarray*}
d\overline \mu _{(\pi (g_1)  , \pi (g_2
  ))} (\beta _{\pi (g_1)} + \beta _{\pi (g_2 )} )  & = &  d\overline
  \mu _{(\pi (g_1 ) , \pi(g_2
  ))} \circ d(\pi \times \pi
)_{g_1 , g_2 )} 
(\gamma _{g_1} + \gamma _{g_2} ) \\ & = & d\pi _{g_1 g_2} \circ d\mu _{(g_1 ,
  g_2 )}  (\gamma _{g_1} + \gamma _{g_2} )   = d\pi _{g_1 g_2} (\gamma
_{g_1 g_2} ) \\ & = & \beta _{\pi (g_1 g_2 )} \; . 
\end{eqnarray*}

\noi Therefore $\beta$ is multiplicative and $(G/K , \beta ) $ is a Poisson
group.  This correspondence $D \mapsto \pi _\star D $ is injective
because $\pi ^\star \pi _\star D = D $.  \\

\noi To complete the proof, it suffices to show that for any Poisson
group structure \\
$L = L(\beta , T^* (G/K ))$  on $G/K$, $(G, \pi^\star L) $ is a Dirac group.  From Lemma~\ref{pdsprop2}, $\pi ^\star L $ is a Dirac structure on $G$  because $L$ is a Dirac
structure on $G/K$. Note that $\pi ^\star L = L(\gamma , U ) $, where $U = Ann(\mathfrak k)$. Since $\mathfrak k$ is $G$-invariant, $U$ is bi-invariant.  Thus, (M1) is satisfied. \\

\noi It remains to show that $\gamma$ is multiplicative.  This follows
in the same way as before. Since $\pi : (G, \pi ^\star L ) \lra (G/K ,
L) $ is Dirac, $d\pi \gamma = \beta$.  Observe that
\begin{eqnarray*}
d\pi _{g_1 g_2} (\gamma _{g_1 g_2} ) & = & 
 \beta _{\pi(g_1 g_2 )} =
  d\overline \mu _{(\pi (g_1 ) ,\pi ( g_2 )
  )} (\beta _{\pi (g_1)} + \beta _{\pi (g_2 )} ) \\
 &  =  & d\overline \mu
  _{(\pi (g_1 ) ,\pi ( g_2 ) )} \circ (\pi \times \pi )_{(g_1 , g_2 )}
  (\gamma _{g_1} + \gamma _{g_2} ) \\ 
 & = & d\pi _{g_1 g_2} \circ d\mu _{(g_1 ,
  g_2 )} (\gamma _{g_1} + \gamma _{g_2} ). 
\end{eqnarray*}
Therefore $\gamma$ is multiplicative.  This proves the first two statements.\\

\noi Now suppose that $G$ is semisimple.  By the first part of this proposition, to prove our correspondence it suffices to show that every Ad-invariant $\mathfrak k \subset \fg$
 is the Lie algebra of some normal closed subgroup.  Any Ad-invariant $\mathfrak k \subset \fg$ is an ideal.  Since $\fg$ is semisimple, $ \fg = \mathfrak k \oplus \mathfrak k ^\perp$,
 where $\mathfrak k ^\perp $ is determined by the Killing form on $\fg$.
  We know that $\mathfrak k ^\perp$ is also an ideal and semisimple.
  The identity component of the centralizer $Z_G(\mathfrak k ^\perp ) $
 is a closed subgroup with Lie algebra $Z_\fg (\mathfrak k ^\perp) = \mathfrak k$.$\square$
\end{pf}

\subsection{Generalized Complex Groups} 



\noi Here we classify generalized complex groups.  In parallel with the work of Gualtieri, we hope that the following result may be helpful in a geometric construction of quantization of holomorphic Poisson groups.

\begin{prop}\label{27jan1} A generalized complex group structure on a group $G$ is equivalent to a holomorphic Poisson group structure on $G$.


\end{prop}

\begin{pf} By Lemma~\ref{24nov3}, $D = L(\beta , U ) $,
  where U is a bi-invariant subbundle of $G \times \cts ^* $.  Then
  $\mathfrak k = Ann(U_e ) $ is an ideal in $\cts$ and $G$-invariant.  Also by Lemma~\ref{24nov3}, $D_e = L(0,U_e)$, so $D_e
  \oplus \overline{D_e} = \cgts $ implies that $\fk \oplus \overline \fk = \cts$.  Therefore $\fk$ gives $G$ the structure of a complex group.  We henceforth identify $TG \simeq G \times \fg$ by left translation.  Let $\JJ : \VV _G \lra \VV _G$ be the map with i-eigenbundle $D$.  Since $\fk \subset D_g$ for any $g$, $\JJ  \fk =  \fk$, and since $\fk \subset \overline D_g$, $\JJ \fk \subset \fk$ because $\overline D_g$ is the $-i$-eigenspace of $\JJ _g$.  Therefore, since $\fk \oplus \overline \fk = \cts$, $\JJ \fg \subset \fg$ and $\JJ $ is of the form
\[
\JJ = 
\left[ 
\begin{array}{cc} 
J & Q \\
0 & -J^* \\
\end{array}
\right]
\]
\noi where $J$ is the complex structure with i-eigenbundle $G \times \fk$.  Gualtieri shows that generalized complex structures of this form are equivalent to a complex structure $J$ 
with holomorphic Poisson structure 
$\rho _\sharp = JQ + iQ$ \cite{gua2} \cite{gua3}.
  One may check that the graph of $ -\tfrac{-i}{2} Q _{|((T^{0,1})^*}$ is contained in the i-eigenbundle of $\JJ$ so that $\be _\sharp  = pr_{|T^{(0,1}}  -\tfrac{-i}{2} Q _{|(T^{0,1})^*} = -\tfrac{-i}{2}(1+iJ) \circ Q _{|(T^{0,1})^*}$, from which one can compute that $\be _\sharp $ is a multiple of $Q + iJQ$.  This implies that $\be$ is multiplicative if and only if $Q$ is multiplicative, whence the desired result follows. $\square$

\end{pf}  

\begin{rem} We have noted that Crainic defines a category of generalized complex manifolds with morphisms which I will call here \emph{Crainic-GC maps}.  When $f$ is a submersion, as is the case for the group multiplication map, $f$ being Crainic-GC implies that $f$ is Dirac, but the converse is not true in general.  However, we do notice the following result.
\end{rem}
\begin{prop} A Crainic-GC group is the same as a generalized complex group. 
\end{prop}
\begin{pf} Vaisman \cite{vai2} shows that groups with generalized complex structure such that multiplication is a Crainic-GC map are equivalent to a complex Poisson groups such that the complex and Poisson structure form a Poisson-Nijenhuis structure. Laurent-Gengoux, Sti\'enon, and Xu \cite{lsx} describe how a Poisson-Nijenhuis structure is equivalent to a holomorphic Poisson structure.  Thus, Vaisman's result implies that a Crainic-GC group structure is equivalent to a holomorphic Poisson group structure.  What Proposition~\ref{27jan1} shows is that Crainic-GC groups are the same as generalized complex groups.$\square$  
\end{pf}



\subsection{Integrability and Dirac Groups}

\noi Here we generalize several formulations of integrability for Poisson groups treated in \cite{luw}, \cite{vai} to the case of almost Dirac groups.  We begin with an independent observation.



\begin{prop}\label{22jan4} Let $(G , D = L(\be, U))$ be an almost Dirac group.  If there is an ideal $V \subset \cts$ complementary to $\fk = Ann(U_e)$ (for example, if $\fg$ is reductive), then $\be : G \lra \wedge ^2 \cts / \fk$ can be 
lifted to $\hat \be : G \wedge ^2 V \subset \wedge ^2 \cts$.  The almost Dirac structure $D$ is integrable if and only if $L(\hat \be , T^*G _\cc)$ is integrable and $d_g \be  (\fk ) = 0$ for all $g \in G$ (i.e. $\LL _{\tilde X} \be = 0$ for all $X \in \fk$, where $\tilde X$ is left-invariant and $\tilde X _e = X$).

\end{prop}

\begin{pf}
First suppose that  $D$ is integrable.  We first show that $\fk (\be ) = 0$, i.e. $d_x \be (\fk ) =0$.  If $\eta$, $\si \in V^* \subset \cts$.  In order for $[ \be _\sharp \eta + \eta , \be _\sharp \si + \si ] $ to be a section of $D$, we need $\phi = [\be _\sharp \eta , \si ] - [\be _\sharp \si , \eta ] $ to be a section of $G \times V^* = Ann (G \times \fk)$.  A simple computation yields that for $Y \in \fk$, $\phi (Y) = -Y(\be(\eta, \si))- \si [\be _\sharp \eta, Y] + \eta [\be _\sharp \si , Y]$.  Another simple computation shows that $-\si[\be _\sharp \eta , Y] = Y (\be(\eta , \si))$ so that $\phi (Y) = Y(\be (\si , \eta))$.  Integrability implies that $Y(\be (\eta , \si)) (g) = d_g \be (Y) (\eta , \si) = 0$.  \\

\noi To check that $D^\prime$ is integrable, it is enough to check this on a frame, and in fact, it is enough to check that $[\hat \be _\sharp \eta + \eta , \al ]$ is a section of $D^\prime $ for constant sections $\al $ of $G \times \fk ^*$ and $\eta $ of $G \times V^*$.  But $[\hat \be _\sharp \eta + \eta , \al ] = [\hat \be _\sharp \eta , \al] = \iota _{\hat \be _\sharp \eta} d \al + \tfrac{1}{2} d\iota_{\hat \be _\sharp \eta} \al = \iota _{\hat \be _\sharp \eta} d \al $.  For a constant section $X$ of $G \times V$, $[\hat \be _\sharp \eta, \al ] (X) = d\al (\hat \be_\sharp \eta , x) = \hat \be _\sharp \eta (\al (X)) - x(\al(\hat \be_\sharp \eta)) + \al ([\hat \be _\sharp \eta, X]) = 0 -0+ \al ( [\hat \be _\sharp \eta , X] ) =0$ because $[\hat \be _\sharp \eta , X] $ is a section of $G\times V$ and $\al \in \fk ^*$.  Therefore $[\hat \be _\sharp \eta + \eta , \al]$ is a section of $G \times \fk ^* \subset D^\prime$.  \\

\noi Now suppose that $D^\prime$ is integrable and $\fk (\be ) =0$.  The first paragraph of this proof shows that for constant sections $\eta $, $\si$ of $G \times V^*$, $[\be _\sharp \eta , \si ] - [\be _\sharp \si , \eta]$ is a section of $G \times V^*$ because $\fk (\be )  = 0$.  Thus, $[\be_\sharp \eta + \eta , \be _\sharp \si + \si ] $ is a section of $D$.  Now to show integrability of $D$, we simply must show that for constant sections $k$ of $G \times \fk$ and $\eta $ of $G \times V^*$, $[k ,\be _\sharp \eta + \eta ]$ is a section of $D$.  Because $\fk$ is an ideal and $\fk (\be) =0$, one may easily check that $[k , \be _\sharp \eta] $ is a section of $G \times \fk$.  Also $[k , \eta](X) =0$ for any constant section $X$ of $G \times \cts$ so that $[k , \eta] = 0$. $\square$

\end{pf} 

\noi When $\fg$ is semisimple, $\fk$ has a unique complementary ideal $\fk ^c$, so $D^\prime$ in Proposition~\ref{22jan4} is canonical.  

\begin{cor} Let $G$ be a semisimple Lie group. The Dirac group structures are parameterized by pairs ($\fk$, $\be$) of a G-invariant ideal $\fk $ of $\cts$ and a complex Poisson group structures $\be : G \lra \wedge ^2 \cts$ on $G$ such that $\be$ takes values in $\wedge ^2 \fk ^c $ and $d_g \be (\fk) = 0$ for all $g \in G$. 
\end{cor}

\begin{lem}\label{22jan3} Given an ideal $\fk \subset \cts$, there is a well defined Schouten bracket 
$[ \be , \be ] \in G \lra  \wedge ^3 (\cts / \fk )^*)$.  In fact, the Schouten bracket $[P,Q]$ makes sense for any $P ,Q : G \lra  \wedge ^2 \cts / \fk )$, as long as $d_g P (\fk ) = d_g Q (\fk) = 0$ for all $g \in G$.
\end{lem}

\begin{pf} 
Let $P$, $Q$ be as above.  By choosing a splitting $\si$ of $\pi : \cts \lra \cts /\fk$, there are sections $\hat P , \hat Q $ of $ G \times \wedge ^2 \si (\cts / \fk ^*) \subset G \times \wedge ^2 \cts)$, so we define the Schouten bracket $[P,Q] = \pi [\hat P , \hat Q]$.  One may verify that since $\fk$ is an ideal, $[P,Q]$ does not depend on the choice of spliting. $\square$
\end{pf}

\noi Let $(G, D = L(\be, U))$ be an almost Dirac group.  Let $\FF \subset \CC ^\infty (G) \otimes \cc $ denote the subsheaf consisting of all functions on which  vectors in $G \times \fk \subset   TG_\cc$ vanish.  Then $\be$ defines a bracket operation $\{ f, g \} = \be (df , dg)$ on $\FF $ with values in $\CC ^\infty (G) _\cc$.  By considering $\hat \be$ as in Lemma~\ref{22jan3}, this bracket can also be extended to one on $\wedge ^2 \CC ^\infty (G) \otimes \cc $.      

\begin{prop}\label{22jan1} Let $(G,D)$ be an almost Dirac group as above. Then the following are equivalent:
\begin{enumerate}
\item $(G, D)$ is a Dirac group (i.e. $D$ is integrable)
\item  The bracket operation $\{ \, , \, \} : \FF \times \FF \lra \CC ^\infty (G)_\cc$ satisfies the Jacobi identity.  In this case, $\{ \; , \}$ is a Lie bracket on $\FF$.  
\item  $[\be , \be ] = 0$.
\end{enumerate}
\end{prop}

\begin{pf} 
(1 $\iff$ 2) First we note that locally there exists a frame of 
$G \times U \subset G \times \cts$
 consisting of sections of the form 
$df$ for $f \in \FF$.  Gualtieri
 shows \cite{gua} that integrability
 of $D$ is equivalent to the Nijenhuis 
operator $Nij (A,B,C) = \tfrac{1}{3}(\langle [A,B],C \rangle + \langle [B,C],A \rangle + \langle [C,A], B \rangle ) $
 vanishing on all sections $A,B,C $ of $D$.  We may, choose sections of the form $X + df$,
 where $f \in \FF (G)$ and $X_{|U} = \be (df, -)$.  If $D$ 
is integrable, then $[X + df , Y + dg] = [X,Y] + d (\be(df,dg))$
 is a section of $D$ so that $\be (df ,dg) \in \FF (G)$.  
The integrability condition is 
$0 = Nij(X + df , Y + dg , Z + dh ) = \{ f , \{ g, h\} \} + \{ g , \{ h,f \} \} + \{ h , \{ f , g \} \}$.  On the other hand, if $\{ \; , \; \}$ satisfies the Jacobi identity, then $Nij_{|D} = 0$, so $D$ is integrable.  Hence $\be (df,dg) \in \FF (G)$ for $f,g \in \FF (G)$, which makes $\FF$ a sheaf of Lie algebras. \\

\noi (2 $\iff$ 3) It follows directly from ~\cite{vai} equation 1.16, originally appearing in \cite{bhv}, that with $\hat \be $ defined as in Lemma~\ref{22jan3}, $f,g,h \in \FF (G)$, $[\be , \be ] (df, dg,dh) = [\hat \be , \hat \be ] (df,dg,dh) = \{ f , \{ g, h\} \} + \{ g , \{ h,f \} \} + \{ h , \{ f , g \} \}$.  $\square$

\end{pf}



\noi The classification \cite{vai} of Poisson group structures in terms of Lie bialgebras can be extended to Dirac group structures.  Proposition~\ref{22jan4} shows that $d_e \be (\fk) = 0$, so $d_e \be: \cts \lra \wedge ^2 \cts / \fk$ factors through $\cts / \fk$ and can be thought of as $d_e \be : \cts / \fk \lra \wedge ^2 \cts \fk$.

\begin{lem}\label{5feb2} Let $G$ be a Lie group, $\fk$ a G-invariant ideal and $\be \in \Ga (G , G \times \wedge ^2 \cts / \fk)$ such that $\be (e) = 0$.  If the almost Dirac structure $L(\be , U = Ann(G \times \fk))$ is integrable, then $d_e \be (\fk ) = 0$, and  $d_e \be ^* : \wedge ^2 \cts / \fk \lra \cts /\fk ^*$ defines a Lie bracket on $\cts  / \fk ^*$.  If $G$ is simply connected and $\be$ is multiplicative, then $L(\be , U = Ann(G \times \fk))$ is integrable if and only if $d_e \be ^*$ is a Lie bracket.
\end{lem}

\begin{pf} 
We only need to make a few modifications to the proof presented in \cite{vai}.  First, as in Proposition~\ref{22jan4}, $d_g \be (\fk) = 0$ is a necessary condition for integrability.  As in \cite{vai}, the bracket $[ \; , \; ] = d_e \be ^* : \wedge ^2 (\cts / \fk ) ^* \lra  (\cts / \fk )^*$ is given by $[\al , \be ] = d(\be(\tilde \al , \tilde \be )) _e$, where $\tilde \al $ and $\tilde \be$ are any sections of $G \times \cts / \fk ^*$ with $\tilde \al (e) = \al $, $ \tilde \be (e) = \be$.  If $D$ is integrable, then by Proposition~\ref{22jan1}, $\{ \; , \; \}$ is a Lie bracket on the sheaf $\FF$.  Let $\FF_e ^i$ denote all germs of functions in the stalk $\FF _e$ which vanish up to order $i$ at $e$.  Then there is an isomorphism $\FF ^0 _e / \FF ^1 _e \tilde \lra U_e \subset \cts ^*$ sending $[f] \mapsto df_e$.  Then $\{f,g\} = \be (df ,dg) \mapsto d(\be(df,dg))_e = [df_e , dg_e]$, so $d_e \be ^*$ satisfies the Jacobi identity.  If, however, we assume that $d_e \be ^*$ satisfies the Jacobi identity, then the proof in \cite{vai} follows as written and $[\be, \be] = 0$, which by Proposition \ref{22jan1} implies that $D$ is integrable.  $\square$

\end{pf}

\begin{thm} If $G$ is simply connected, a Dirac group structure on $G$ is equivalent to an ideal $\fk$ of $\cts$ and a Lie bialgebra structure on $\cts / \fk$.  
\end{thm}

\begin{pf} First let $L(\be , U)$ be a Dirac group structure on $G$. For the ideal $\fk = Ann (U_e)$, a Lie bi-algebra structure on $\cts / \fk$ is equivalent to a cocylce $\epsilon : \cts / \fk \lra \wedge ^2 \cts / \fk$ such that $\epsilon ^*$ is a Lie bracket.   Then by Proposition~\ref{5feb1}, $d_e \be : \fg \lra \wedge ^2 \cts / \fk$ is a cocycle.  Since $L(\be , U)$ is integrable, $d_e \be (\fk) = 0$ by Lemma~\ref{22jan4}.  So $d_e \be$ factors to some cocycle $\epsilon : \cts / \fk \lra \wedge ^2 \cts / \fk$.  By Lemma \ref{5feb2}, integrability of $L(\be , U)$ implies that $\vep ^*$ satisfies the Jacobi identity.  \\

\noi Now let $\epsilon : \cts / \fk \lra \wedge ^2 \cts / \fk$ define a Lie bi-algebra structure.  Then letting $\epsilon ^\prime  : \fg \lra \wedge ^2 \cts / \fk$ be the composition of $\epsilon $ with $\fg \hookrightarrow \cts \lra \cts / \fk$.  Then $\epsilon ^\prime$ is a cocycle for $\fg$.  By Proposition \ref{5feb1}, there is a unique $\be : G \lra \wedge ^2 \cts / \fk$ such that $d _e \be = \epsilon ^\prime$.  Since $\epsilon$ satisfies the Jacobi identity, $L(\be , Ann(\fk))$ is integrable (by Lemma \ref{5feb2}).  \\

\noi Clearly these two processes are inverses of each other. $\square$  

\end{pf}

\subsection{Dual-Dirac Groups}

\noi A \emph{dual-Dirac group} is a group with Dirac structure such that group multiplication is a dual-Dirac map.

\begin{prop} Dual-Dirac group structures on a Lie group $G$ are in bijection with pairs $(E, \vep)$ of a G-invariant ideal $E $ of $\cts$ and $\vep \in \wedge ^2 E^*$ such that $\vep$ is G-invariant and $\vep $ is a cocycle in the Lie algebra cohomology of $E$ with values in $\cc$ (i.e. $\vep (x, [y,z] ) + \vep (y, [z,x]) + \vep (z , [x,y]) = 0 $ for all $x,y,z \in E$).   If $G$ is connected semisimple, dual-Dirac group structures on $G$ are given by ideals of $\cts$.  
\end{prop}

\begin{pf}  Let $(G, D)$ be an almost dual-Dirac group.  That is, multiplication is dual-Dirac and $D$ is an almost Dirac structure.  Then poinwise, $D_g = L(E_g , \vep _g )$.  We again identify $TG _\cc \simeq G \times \cts$ by left-translation, and we think of $E $ as an assignment of a subspace of $\cts$ for each $g \in G$.  Condition $(M1)^*$ is that $Ad(h\inv)E_g + E_h \subset E_{gh}$.  Hence, $E_h \subset E_{gh}$ for all $g$, $h$.  Clearly then $E_g = E_e$ for all $g$.  Thus, $E$ is constant or left-invariant.  But also $Ad(h\inv ) E \subset E$ for all $h \in G$, so $E$ is a G-invariant ideal of $\cts$.  The requirement $(M2)^*$ is that $(Ad(h\inv)^* \vep _{gh} , \vep _{gh} ) = (\vep _g , \vep _h )$, so $\vep : G \lra \wedge ^2 E^*$ is also constant.  Now $(M2 ^*) $ holds if and only if $Ad (G)^* \vep = \vep$.  \\

\noi Gualtieri shows that almost Dirac structures of the form $L(E, \vep)$ are integrable if and only if $E$ is an integrable distribution and $d_E \vep = 0$.  Since $E$ is an ideal, $G \times E \subset G \times \cts \simeq TG _\cc$ is an integrable distribution.  \\

\noi When $G$ is connected semisimple, $E$ is also semisimple, so $H^2 (E , \cc) = 0$ \cite{wei}.  Hence, $\vep = \phi \circ [ \,  , \, ]$ for some $\phi \in \cts ^*$.  $G$-invariance of $\vep$ implies that $\phi$ is $G$-invariant.  Since $G$ is semisimple, this is only possible if $\phi = 0$.  $\square$

\end{pf}

\section{Twisted Groups}

\noi Here we define groups in $\TT \DD ^*$ and $\TT \DD$ and describe the data for such groups.  

\subsection{Groups in $\TT \DD ^*$ }
\noi We consider only Courant algebroids of the form $\VV _{G , H}$, i.e. exact Courant algebroids in the category $\TT \DD ^* \subset \CC \DD ^*$.  
Products do not exist in $\CC$ or $\CC \DD ^*$, so we cannot have groups objects in these categories.  However, we will say that a \emph{group} in $\TT \DD ^*$ consists of a group $G$, a Dirac structure $D \subset \VV _{G, H}$, and a 2-form $B \in \Om ^2 (G\times G)$ such that $(\mu ,B)  : (G \times G , D \oplus D ) \lra (G , D)$ is a morphism in $\TT \DD ^*$.  In Proposition~\ref{4augb} we trivialize the vector bundle $\VV _{G, H}$ by left translation and let $\pi _1 \; , \pi _2 : G \times G \lra G$ be the first and second projection maps respectively.

\begin{prop}\label{4augb} A group $(G,D,H,B)$ in $\TT \DD ^*$ is given by the following data: 
\begin{enumerate}
\item a G-invariant ideal $E \subset \cts$,
\item a section $\vep \in \Ga (G, \wedge ^2 \tilde E)$, where $\tilde E$ is the bi-invariant distribution defined by $E$,
\item a section $B = B_1 + B_2 + 0 $ of $\wedge ^2 T^*G =( \pi _1 ^* \wedge ^2 T^*G )\oplus (\pi _2 ^* \wedge ^2 T^*G )\oplus (\pi _1 ^* T^*G \otimes \pi _2 ^* T^*G )$ such that on $\tilde E$ one has $(B_1)_{(g,h)}  = \vep _g - Ad (h\inv )^* \vep _{gh} $, $(B_2)_{(g,h)} = \vep _h - \vep _{gh}$ for all $g,h \in G$.  Similarly, we also require that $(dB_1)_{(g,h)}  = H _g - Ad (h\inv )^* H _{gh} $, $(dB_2)_{(g,h)} = H _h - H _{gh}$ for all $g,h \in G$, 
\item $H_{|E} = -d_E \vep$, where $d_E \vep$ is the restriction to $\tilde E$ of the differential of any 2-form restricting to $\vep$.  
\end{enumerate}
\noi For such a $(E, B, \vep , H)$, $D = L(\tilde E , \vep) \subset \VV _{G,H}$ is the corresponding Dirac structure on $G$.  If we require that  $B = \pi _1 ^* \om \oplus \pi _2 ^* \om $ for some $\om \in \Om ^2 (G)$, then $\om _{|E} = 0$ and $\vep $ is constant (i.e. bi-invariant).  In particular, if $G$ is semisimple, then $\vep = 0$.


\end{prop}
\begin{pf} We use left translation to trivialize $\VV _G$, and pointwise, there is $E_g \subset \cts$ and $\vep _g \in \wedge ^2 E_g ^*$ such that $D_g = L(E_g , \vep _g)$.  The condition $(M1)^*$ is that $Ad(h)\in E_g + E_h \subset E_{gh}$, which implies that $E := E_e = E_g $ for all $g \in G$ and $Ad(h \inv )E \subset E$.  Therefore, $E$ is a $G$-invariant ideal.  Condition $(M2)^*$ is that $ (Ad(h\inv)^*\vep _{gh} , \vep _{gh} ) = (\vep _g , \vep _h ) + B_{(g,h)}$, which is exactly condition (3) of this proposition.  Therefore, for $(\mu , B)$ to be a morphism in $\TT \DD ^*$, it is necessary and sufficient that $\mu ^* H = H \oplus H + dB$ as is explained in part (3) and that $E$ ,$\vep$, and $B$ be as in (1) and (2).  For $D = L(E, \vep) \subset \VV _{G, H}$ to be integrable, $d_E \vep = -H_{|E}$.  \\

\noi Now if $B = \om + \om$, then (2) implies that $\om _{|E} = 0$, $\vep$ is constant and $G$-invariant.  When $\cts$ is semisimple, so is $E$.  The only possible $\vep$ is $\vep = 0$.  $\square$

\end{pf}

\noi There is an example of a group in $\TT \DD ^*$, called the Cartan-Dirac structure, described in detail in \cite{abm}.

\subsection{Groups in $\TT \DD$}
\noi Here we define twisted Dirac groups and show that they are a generalization of twisted Poisson groups.  

\begin{Def} We define the full subcategory $\TT \DD$ of $\CC \DD$  of all objects  of the form  $(M , \VV _{M,H} , id , D) $.  The data for an object in $\TT \DD$ is just a triple $(M,H,D)$.  We say that $(G,H,D)$ is a \emph{group in $\TT \DD$} if $G$ is a Lie group with multiplication $\mu$ such that $(\mu , B) : (G\times G , H \oplus H , D \oplus D ) \lra (G, H,D)$ is a map in $\TT \DD$ for some 2-form $B$ on $G \times G$.   
\end{Def}

\begin{prop}  The object $(G , D \subset \VV _{M, H})$ in $\TT \DD$ is a group in $\TT \DD$ if and only if $(G,D)$ is an almost Dirac group and $H $ is a section of $\wedge ^2 U \subset \wedge ^2 T^*G$ such that $[\be , \be ] = \be _\sharp H$ and such that $\mu ^* (H \oplus H) -H$ is exact.   
\end{prop}

\begin{pf} The condition for $\mu $ to be a map in $\TT \DD$ is simply for $\mu$ to be a Dirac map.  This happens exactly when $(G,D)$ is an almost Dirac group.  In this case, $D = L(\be , U)$.  Integrability of $D \subset \VV _{M,H}$ is equivalent to
 $Nij_H = 0$ on $D$, where $Nij_H (X + \xi , Y + \eta , Z + \al ) = Nij(X + \xi , Y + \eta , Z + \al ) + H(X,Y,Z)$.  We recall \cite{gua} that $Nij_H (X + df ,Y + dg ,Z + dh) = Jac_{\{\; , \; \}} (f,g,h) + H(X, Y,Z)$.  What this says is 
that when
 $X _{|U} = \be _\sharp (df)$, $Y _{|U} = \be _\sharp (dg)$,
 and $Z _{|U} = \be _\sharp (dh)$, $H(X,Y,Z) = -Jac_{\{\; , \; \}} (f,g,h) $.
  This means that $H(X,Y,Z)$
 really only depends on the images 
$\overline X $, $\overline Y$,
 and $\overline Z$ in $U^*$,
 not on $X$, $Y$, $Z$ themselves.
  In other words, $H \in \wedge ^2 U \subset \wedge ^2 T^*G$.
  Now, $H(X,Y,Z) = H(\be _\sharp (df) , \be _\sharp (dg) , \be _\sharp (dh)) = \be _\sharp H (df , dg, dh)$. 
 The integrability condition is that $\be _\sharp H = Jac _{\{ \, , \, \}} = [\be , \be ]$. The condition on $B$ and $H$ being related by the pullback $\mu ^*$ is unrelated to either of these conditions.  $\square$
\end{pf}

\end{document}